\documentclass[11pt]{amsart}
\usepackage[a4paper,margin=2.6cm]{geometry}
\usepackage{graphicx}
\usepackage{amssymb,amsmath}
\usepackage{booktabs}
\usepackage{array}
\usepackage{longtable}
\usepackage{enumitem}
\usepackage{xcolor}
\usepackage[colorlinks=true,linkcolor=blue!50!black,citecolor=blue!50!black,urlcolor=blue!50!black]{hyperref}

\newtheorem{theorem}{Theorem}[section]
\newtheorem{proposition}[theorem]{Proposition}
\newtheorem{lemma}[theorem]{Lemma}

\theoremstyle{definition}

\theoremstyle{remark}
\newtheorem{remark}[theorem]{Remark}

\newcommand{\D}{\mathbb{D}}
\newcommand{\R}{\mathbb{R}}

\newcommand{\Q}{\mathbb{Q}}
\newcommand{\Z}{\mathbb{Z}}
\newcommand{\PSL}{\mathrm{PSL}}
\newcommand{\PGL}{\mathrm{PGL}}
\newcommand{\Aut}{\mathrm{Aut}}
\newcommand{\Isom}{\mathrm{Isom}}
\newcommand{\Dl}{\Delta}
\newcommand{\Dpm}{\Delta^{\pm}}
\newcommand{\Dist}{\mathcal{D}}
\newcommand{\rc}{r_c}

\title[Extremal disc packings in surfaces of minimal genus]{Extremal disc packings in compact hyperbolic surfaces of minimal genus}
\author{Ernesto Girondo}
\address{Departamento de Matem\'aticas, Universidad Aut\'onoma de Madrid and Instituto de Ciencias Matem\'aticas ICMAT (CSIC-UAM-UCM-UC3M), 28049 Madrid, Spain}
\email{ernesto.girondo@uam.es}

\author{Cristian Reyes-Monsalve}
\address{Departamento de Ciencias Exactas, Facultad de Ingenier{\'\i}a, Universidad San Sebasti\'an, Concepci\'on, Chile}
\email{cristian.reyes@uss.cl}
\begin{document}


\begin{abstract}
A compact hyperbolic surface is $k$-extremal if it contains $k$ disjoint embedded discs of the largest possible radius; such surfaces are uniformized by surface subgroups of the triangle groups $\Dpm(2,3,N)$. We determine all extremal $k$-packings in every $k$-extremal surface of minimal genus, orientable of genus $2$ and non-orientable of genus $3$: $373$ surface groups in total, all treated in a unified way; $264$ of them had never been studied before. The computation is automatic and verified in exact arithmetic, each packing being certified by an element $\alpha\in\PGL_2$ over an explicit number field with $\alpha^{-1}K\alpha<\Dpm(2,3,N)$. We also obtain the full isometry groups of these surfaces. In every surface but one, all the extremal packings are related by isometries. The exception is the genus two surface $X_6$ of Girondo and Gonz\'alez-Diez, which has exactly five extremal $4$-packings, one invariant under all isometries and four forming an orbit; it is also the only one of these surfaces uniformized by two non-conjugate subgroups of the same triangle group, and we show that the two facts are equivalent. We find its equation, $y^2=x^6+2\sqrt{3}\,i\,x^3+1$. The results form a catalogue with detailed data and pictures of every surface (packings, isometries, Weierstrass points, Dirichlet domains), distributed with programs that allow the reader to consult it and to reproduce it completely.
\end{abstract}

\maketitle

\section{Introduction}

A packing of metric discs in the hyperbolic plane has the highest possible density when every point outside the discs lies in a region bounded by three pairwise tangent discs of the packing (see \cite{Boroczky_1978}). These optimal \emph{B\"or\"oczky packings} are closely related to the \emph{extremal surfaces}: a compact hyperbolic surface $X$ is extremal if it contains an embedded disc that lifts to a B\"or\"oczky packing of the universal cover. Such \emph{extremal discs} realize the maximum radius of an embedded disc in any surface with the same topology, which is determined by the Euler characteristic $\chi(X)$ and the orientability character. Extremal surfaces exist in every genus, both in the orientable and in the non-orientable case; see \cite{Girondo_Gonzalez-Diez_1999, Girondo_Gonzalez-Diez_2002, Girondo_Gonzalez-Diez_2002_2, Girondo_Nakamura_2007} for an overview.

The theory extends naturally to collections of $k>1$ disjoint isometric embedded discs, the so-called $k$\emph{-packings}. It turns out that for certain triples $(\chi,k,\pm)$ no $k$-packing of a compact hyperbolic surface with Euler characteristic $\chi$ and the given orientability character can lift to a B\"or\"oczky packing, and nothing is known about the extremal radius in those cases. The necessary and sufficient condition for the existence of \emph{extremal $k$-packings} is that $k$ divides $6\chi$ (see \cite{Girondo_2018, deBlois_2020, Girondo_Reyes_2021, Girondo_Reyes_2023}).

The aim of this paper is to extend the database of classical extremal surfaces (completed in the series of papers \cite{Girondo_Gonzalez-Diez_2002_2, Girondo_Nakamura_2007, Nakamura_2005, Nakamura_2009, Nakamura_2012, Nakamura_2013, Nakamura_2016}) to $k$-extremal surfaces of the lowest possible genus for every admissible $k$, that is, to orientable surfaces of genus $2$ and to non-orientable surfaces of genus $3$. The admissible values of $k$ are $k\in\{1,2,3,4,6,12\}$ in the orientable case and $k\in\{1,2,3,6\}$ in the non-orientable one. Since the orientable double cover of a $k$-extremal non-orientable surface of genus $g$ is a $2k$-extremal orientable surface of genus $g-1$ (see \cite{Girondo_Reyes_2023}), it is natural to treat both families at the same time.

A $k$-extremal surface is uniformized by a subgroup $K$ of a triangle group $\Dpm(2,3,N)$, so it comes with an obvious extremal $k$-packing: $K$ preserves the tessellation of the hyperbolic plane by regular $N$-gons associated with the triangle group, which therefore induces a decomposition of the surface into $k$ regular $N$-gons, and the discs inscribed in them form the packing. For each surface we determine the \emph{complete} list of extremal packings, not only this one. The main difficulty is that all the triangle groups involved are arithmetic, so a surface may contain \emph{hidden} packings, whose lifts to the hyperbolic plane are not related to the original tessellation by an element of the triangle group. In \cite{Girondo_Reyes_2023} and in \cite{Reyes_thesis} such hidden packings were detected by a combination of numerical experiments and a final algebraic verification, carried out surface by surface. Chapter~7 of \cite{Reyes_thesis} describes a complete procedure of this kind (a list of admissible distances, the candidates inside the circumscribed disc of the central polygon, a numerical test, the construction of an automorphism of the surface and an algebraic verification), and applies it to the $50$ non-orientable surfaces of genus $3$ and to their double covers. Here we turn that procedure into an algorithm (Section~\ref{sec:algorithm}) whose output is verified in exact arithmetic, and we run it on all the $373$ groups involved. For $109$ of them the extremal packings were already known: the $50$ non-orientable groups and their $50$ canonical double covers \cite{Reyes_thesis} (the $11$ with $k=1$ also in \cite{Girondo_Nakamura_2007}), and the $9$ groups with $k=1$ in genus $2$ \cite{Girondo_Gonzalez-Diez_2002_2, Nakamura_2002}; the classification of the other $264$ is new. We treat all of them in a unified way, with the same algorithm and the same exact verification, so that the results can be compared across the whole list: this is what allows us to locate the double covers among the orientable surfaces, to detect the surfaces that are $k$-extremal for two values of $k$, and to identify the groups that uniformize the same surface. Besides the number and position of all the packings, the method produces, from the same exact data, the full automorphism group of each surface (including the automorphisms that are not induced by the triangle group), decides whether a surface is real, and detects which uniformizing groups are conjugate in $\PSL_2(\R)$ without being conjugate in the triangle group.

Our main results are collected in Section~\ref{sec:results}. In the non-orientable case every $k$-extremal surface of genus $3$ admits either one or two extremal $k$-packings, and the orientable double cover always has exactly as many extremal $2k$-packings (Theorems~\ref{thm:nonor} and~\ref{thm:cover}; this is special to genus~$3$, as shown by five non-orientable surfaces of genus~$4$, Remark~\ref{rem:genus4}); this confirms by an independent method the results obtained in \cite{Reyes_thesis} and announced in \cite{Girondo_Reyes_2023}. In the orientable case the number of packings is again one or two except for two surfaces of genus $2$ (Theorem~\ref{thm:or}), one of them, the surface $X_6$ of \cite{Girondo_Gonzalez-Diez_2002_2}, containing four extremal discs when regarded as a $1$-extremal surface and five extremal $4$-packings when regarded as a $4$-extremal surface. A natural question, implicit in the method of \cite{Girondo_Reyes_2023} where hidden packings were always found as images of the original one under an automorphism, is whether two extremal packings of the same surface are always related by an isometry. The answer is affirmative for all the surfaces considered here with the single exception of $X_6$. That this surface carries two extremal $4$-packings not related by any automorphism was already observed in Section 4.5 of \cite{Girondo_2018} (the discs are inscribed in the faces of the two non-isomorphic uniform dessins of type $(2,3,9)$ found in \cite{Girondo_Torres_Wolfart_2012}); here we obtain the complete picture: $X_6$ has exactly five extremal $4$-packings, which split into an orbit of four and an isolated packing invariant under the whole isometry group (Theorem~\ref{thm:or}(vi)). The surface $X_6$ is also the only one in our list for which two different conjugacy classes of groups in a triangle group give rise to the same Riemann surface, and the two facts are equivalent: the orbits of extremal packings under the isometry group correspond to the conjugacy classes, in the triangle group, of the groups uniformizing the surface (Lemma~\ref{lem:orbits}). We identify $X_6$ as the genus two curve
$$y^2=x^6+2\sqrt{3}\,i\,x^3+1\ \ \cong\ \ y^2=(x^3-1)(x^3+7+4\sqrt{3}),$$
whose group of conformal automorphisms is the dihedral group of order $12$ (Theorem~\ref{thm:X6}). It was known from \cite{Girondo_Gonzalez-Diez_2002_2} to have exactly four extremal discs, permuted transitively by its automorphism group $D_6$, to be a curve $y^2=(x^3-1)(x^3-\lambda^6)$ for some $\lambda$, and its six Weierstrass points had been located there with respect to the extremal disc; our results determine $\lambda^6=-(7+4\sqrt3)$ and describe the Weierstrass points in terms of the tessellation of the invariant $4$-packing.

The result is a catalogue of the $423$ surfaces involved (the $373$ classes and the $50$ canonical double covers), described in Section~\ref{sec:data}. For each surface it contains
\begin{itemize}[leftmargin=2em]
\item all its extremal packings, each one with the exact algebraic identity that certifies it and the position of the centres of its discs;
\item its full isometry group and its group of conformal automorphisms, with their isomorphism types, including the automorphisms not induced by the triangle group, and whether the surface is real;
\item the fixed points and mirrors of all its isometries; in particular, for each of the $323$ classes of genus $2$, the position of its six Weierstrass points with respect to its extremal discs;
\item for each extremal packing, the Dirichlet domain of the surface centred at one of its discs, with its side pairings;
\item the complete record of the search of Section~\ref{sec:algorithm}: the admissible distances, every intersection point of the ED-cycles examined, and the reason why each one is kept or discarded;
\item its relations with the other surfaces of the list: canonical double covers, classes uniformizing the same surface, surfaces that are $k$-extremal for two values of $k$.
\end{itemize}
The catalogue is distributed with the programs that produce it. They serve two purposes. A reader can consult all this information without running anything (tables, text files and pictures of every surface, together with a program that identifies in the catalogue any surface given by its monodromy); and a reader can check that all of it is correct: the exact identities certifying every packing can be verified directly, and everything can be produced again, from the enumeration of the groups to the last picture, and compared with the distributed data file by file, with step-by-step instructions that assume no knowledge of the programs.

\section{Extremal packings and triangle groups}\label{sec:groups}

\subsection{The setting}
Let $X=\D/K$ be a compact hyperbolic surface, orientable or not, with Euler characteristic $\chi<0$, uniformized by a torsion-free discrete group $K$ of isometries of the unit disc $\D$ (a Fuchsian group if $X$ is orientable, a proper NEC group otherwise). An extremal $k$-packing in $X$ lifts to a B\"or\"oczky packing of $\D$ whose Dirichlet--Voronoi cells form the regular tessellation $\{N,3\}$ of $\D$ by regular $N$-gons with angles $2\pi/3$, where
$$N=6-\frac{6\chi}{k}.$$
The group $K$ preserves this tessellation, that is, $K$ is contained in the extended triangle group $\Dpm=\Dpm(2,3,N)$, the full group of isometries of the tessellation, with index $12(k-\chi)$. When $X$ is orientable, $K$ is contained in the Fuchsian triangle group $\Dl=\Dl(2,3,N)$, with half the index. The fundamental domain of $K$ can be taken to be a union of $k$ polygons of the tessellation, one for each disc of the packing. We refer to \cite{Girondo_2018, Girondo_Reyes_2023} for the details.

The cases relevant for this paper are listed in Table~\ref{ta:cases}. We write $(\pm,k,N,d)$ with $d$ the index of $K$ in $\Dl$ in the orientable case, $d=[\Dl:K]=6(k-\chi)$, and in $\Dpm$ in the non-orientable case, $d=[\Dpm:K]=12(k-\chi)$. With this convention a non-orientable surface and its canonical double cover (uniformized by $K\cap\Dl$) have the same $N$ and the same $d$: the double covers of the surfaces of a non-orientable case $(-,k,N,d)$ of Table~\ref{ta:cases} belong to the orientable case $(+,2k,N,d)$.

\begin{table}[h]
\centering
\begin{tabular}{cccc@{\qquad\qquad}cccc}
\toprule
\multicolumn{4}{c}{orientable, genus $2$ ($\chi=-2$)} & \multicolumn{4}{c}{non-orientable, genus $3$ ($\chi=-1$)}\\
$k$ & $N$ & $d=[\Dl:K]$ & polygons & $k$ & $N$ & $d=[\Dpm:K]$ & polygons\\
\midrule
$1$ & $18$ & $18$ & one $18$-gon & $1$ & $12$ & $24$ & one $12$-gon\\
$2$ & $12$ & $24$ & two $12$-gons & $2$ & $9$ & $36$ & two $9$-gons\\
$3$ & $10$ & $30$ & three $10$-gons & $3$ & $8$ & $48$ & three $8$-gons\\
$4$ & $9$ & $36$ & four $9$-gons & $6$ & $7$ & $84$ & six $7$-gons\\
$6$ & $8$ & $48$ & six $8$-gons & & & & \\
$12$ & $7$ & $84$ & twelve $7$-gons & & & & \\
\bottomrule
\end{tabular}
\medskip
\caption{The ten cases covered by this paper. The non-orientable cases had been classified in \cite{Reyes_thesis} (and the case $k=1$ in \cite{Girondo_Nakamura_2007}), and the orientable case $k=1$ in \cite{Girondo_Gonzalez-Diez_2002_2, Nakamura_2002}; in these cases the algorithm recovers the known results (Theorem~\ref{thm:nonor} and Remark~\ref{rem:ggd}). In the orientable cases $k=2,3,4,6,12$, the $50$ canonical double covers of the non-orientable surfaces had been classified in \cite{Reyes_thesis}, and some further examples with two packings had been given in \cite{Girondo_Reyes_2023, Reyes_thesis}; the classification of the remaining surfaces is new.}\label{ta:cases}
\end{table}

Throughout the paper we use the following presentation of the extended triangle group,
$$\Dpm(2,3,N)=\langle a,b,c \mid a^2=b^2=c^2=(ab)^2=(bc)^3=(ca)^N=1\rangle,$$
where $a,b,c$ are the reflections in the sides of a triangle $T_0$ with angles $\pi/N$ (at the centre $O$ of a polygon of the tessellation), $\pi/2$ (at the midpoint $C$ of an edge) and $\pi/3$ (at a vertex $A$): $a$ is the reflection in the side $OC$, $b$ is the reflection in the edge containing $C$ and $c$ is the reflection in the side $OA$. Thus $ca$ is the rotation by $2\pi/N$ about $O$, $ab$ is the half-turn about $C$ and $bc$ is the rotation by $2\pi/3$ about $A$. The Fuchsian triangle group is $\Dl=\langle x,y\mid x^2=y^3=(xy)^N=1\rangle$ with $x=ab$, $y=bc$.

Each triangle of the tessellation $\Dpm(T_0)$ of $\D$ by copies of $T_0$ is $g(T_0)$ for a unique $g\in\Dpm$, and the triangles modulo $K$ correspond to the right cosets $Kg$. The permutation action of $a,b,c$ on $\Dpm/K$ is the only piece of information about $K$ that our algorithm needs, and it is also the form in which GAP \cite{GAP4} returns the subgroups.

\subsection{The groups and their conjugacy classes}\label{sec:classes}
Two groups $K,K'<\Dpm$ that are conjugate in $\Dpm$ uniformize isometric surfaces in which the packings correspond; the converse is not true (this is the subject of Section~\ref{sec:hidden}). Table~\ref{ta:all_info} gives the numbers of conjugacy classes of torsion-free subgroups of the relevant index, computed with GAP by the low index subgroups algorithm with the torsion elements $a,b,c,ab,bc,(ca)^j$ (respectively $x,y,(xy)^j$) excluded. On a standard computer each of these enumerations takes at most about a minute, the longest one being the enumeration of the $155$ classes of index $84$ in $\Dl(2,3,7)$.

In the orientable case we have to distinguish between conjugacy in $\Dl$ and conjugacy in $\Dpm$. Let $\{K_1,\dots,K_n\}$ be representatives of the $\Dpm$-conjugacy classes of torsion-free subgroups of index $d$ of $\Dl$, and let $\sigma\in\Dpm\smallsetminus\Dl$. If the normalizer of $K_i$ in $\Dpm$ is contained in $\Dl$, that is, if the surface $\D/K_i$ admits no anticonformal automorphism induced by the triangle group, then $K_i$ and $\sigma K_i\sigma^{-1}$ are not conjugate in $\Dl$, and they uniformize mirror images of each other; otherwise they are $\Dl$-conjugate. It follows that
$$\#\{\Dl\text{-classes}\}=n+r,$$
where $r$ is the number of classes without an induced anticonformal automorphism. In Table~\ref{ta:all_info} the column ``classes in $\Dpm$'' gives $n$, the column ``with anticonf.\ aut.'' gives $n-r$, the number of classes \emph{with} an induced anticonformal automorphism, and the column ``classes in $\Dl$'' gives $n+r$. The numbers $n-r$ and $n+r$ can be read off from the permutation representation: the automorphisms of $\D/K$ induced by $\Dpm$ correspond to the cosets $Kg$ whose stabilizer is $K$ itself, and $g$ reverses the orientation if and only if $g$ is represented by a word of odd length.

We refer to the classes as in the data files of Section~\ref{sec:data}: \texttt{fuchs\_N9\#21}, for instance, is the $21$st class of torsion-free subgroups of $\Dl(2,3,9)$, up to conjugacy in $\Dl$ (an orientable surface), in the order in which GAP lists them, and \texttt{nec\_N9\#2} is the second class of torsion-free NEC subgroups of $\Dpm(2,3,9)$, up to conjugacy in $\Dpm$ (a non-orientable surface).

\begin{table}[h]
\centering
\begin{tabular}{cccccccc}
\toprule
$\Dpm(2,3,N)$ & $k$ & index & classes in $\Dl$ & classes in $\Dpm$ & with anticonf.\ aut. & $k$ & NEC classes\\
\midrule
$(2,3,7)$ & $12$ & $84$ & $155$ & $90$ & $25$ & $6$ & $12$\\
$(2,3,8)$ & $6$ & $48$ & $77$ & $57$ & $37$ & $3$ & $23$\\
$(2,3,9)$ & $4$ & $36$ & $37$ & $24$ & $11$ & $2$ & $4$\\
$(2,3,10)$ & $3$ & $30$ & $20$ & $15$ & $10$ & -- & --\\
$(2,3,12)$ & $2$ & $24$ & $25$ & $20$ & $15$ & $1$ & $11$\\
$(2,3,18)$ & $1$ & $18$ & $9$ & $8$ & $7$ & -- & --\\
\bottomrule
\end{tabular}
\medskip
\caption{Conjugacy classes of torsion-free subgroups. The first block refers to Fuchsian subgroups of $\Dl(2,3,N)$ of the given index (orientable surfaces of genus 2) and the last two columns to proper NEC subgroups of $\Dpm(2,3,N)$ of the same index, now taken in $\Dpm(2,3,N)$ (non-orientable surfaces of genus 3, with half the value of $k$).}\label{ta:all_info}
\end{table}

\section{Two observations about B\"or\"oczky packings}\label{sec:observations}

Every extremal packing in $X=\D/K$ is determined by an inclusion $K<\alpha\Dpm\alpha^{-1}$ for some isometry $\alpha$ of $\D$: the centres of the discs of the lifted packing are the points $\alpha g(O)$, $g\in\Dpm$. Such an $\alpha$ always belongs to the commensurator of $\Dpm$: the group $\Dpm\cap\alpha\Dpm\alpha^{-1}$ contains $K$, and therefore has finite index both in $\Dpm$ and in $\alpha\Dpm\alpha^{-1}$. For a non-arithmetic triangle group the commensurator is $\Dpm$ itself, so that $\alpha\Dpm\alpha^{-1}=\Dpm$ and the packing is unique. But \emph{all} the groups in Table~\ref{ta:cases} are arithmetic, and hidden packings do occur, as was pointed out in \cite{Girondo_2018} for Klein's quartic and for the surface $X_6$ of \cite{Girondo_Gonzalez-Diez_2002_2} (see also Figures 5 and 7 of \cite{Girondo_Torres_Wolfart_2012}).

\subsection{Admissible distances}
The distance between two centres of the lifted packing is the distance between two centres of the tessellation $\{N,3\}$, so it belongs to the discrete set
$$\Dist_N=\{\,\rho(O,g(O)) \mid g\in\Dpm(2,3,N)\,\}\smallsetminus\{0\}=\{d_0<d_1<d_2<\cdots\},\qquad d_0=2R,$$
where $R$ is the inradius of the polygon ($\cosh R=1/(2\sin(\pi/N))$), $\rc$ is its circumradius ($\cosh\rc=\cot\frac{\pi}{N}\cot\frac{\pi}{3}$; when several values of $N$ are involved we write $R_N$ and $\rc(N)$) and $\rho$ is the hyperbolic distance. In particular, if $z$ is the centre of a disc of an extremal packing of $X$, then $\rho(z,\kappa z)\in\Dist_N$ for every $\kappa\in K$ (\cite[Lemma~4.1]{Reyes_thesis}, where the elements of $\Dist_N$ are called \emph{admissible distances}). Decimal values of the first elements of $\Dist_N$ are given in Table~\ref{ta:dist}.

\begin{table}[h]
\centering\small
\begin{tabular}{ccc@{\qquad}cccccc}
\toprule
$N$ & $R$ & $\rc$ & $d_0$ & $d_1$ & $d_2$ & $d_3$ & $d_4$ & $d_5$\\
\midrule
$7$ & $0.545275$ & $0.620672$ & $1.09055$ & $1.80760$ & $2.14081$ & $2.72010$ & $2.89815$ & $3.15551$\\
$8$ & $0.764285$ & $0.860706$ & $1.52857$ & $2.44845$ & $2.91404$ & $3.05714$ & $3.63158$ & $3.94035$\\
$9$ & $0.927539$ & $1.035886$ & $1.85508$ & $2.89096$ & $3.43821$ & $3.68102$ & $4.23209$ & $4.62526$\\
$10$ & $1.061279$ & $1.177334$ & $2.12255$ & $3.23384$ & $3.83599$ & $4.14771$ & $4.24510$ & $4.67894$\\
$12$ & $1.276687$ & $1.401994$ & $2.55337$ & $3.75563$ & $4.42553$ & $4.82307$ & $5.03827$ & $5.10675$\\
$18$ & $1.719107$ & $1.855077$ & $3.43821$ & $4.74604$ & $5.49629$ & $5.99547$ & $6.34485$ & $6.46676$\\
\bottomrule
\end{tabular}
\medskip
\caption{Inradius $R$ and circumradius $\rc$ of the polygons of the tessellation $\{N,3\}$, and the first six admissible distances $d_0<\dots<d_5$, rounded to six decimals. The exact values are given by $\cosh R=1/(2\sin\frac{\pi}{N})$, $\cosh\rc=\cot\frac{\pi}{N}\cot\frac{\pi}{3}$ and, for the admissible distances, by the elements $\tanh^2\frac{d_i}{2}\in L_N$ computed in Step~2 of Section~\ref{sec:algorithm}.}\label{ta:dist}
\end{table}

The algorithm needs the elements of $\Dist_N$ up to a certain bound $D$: a list $d_0<d_1<\dots<d_m\le D$ together with the guarantee that no other element of $\Dist_N$ is at most $D$. Obtaining such a complete list is not completely trivial, because the distance from $O$ to the centres of the polygons is not monotone with respect to the combinatorial distance in the tessellation. This difficulty was pointed out in Chapter~7 of \cite{Reyes_thesis}, where it is shown that the distances to the polygons at combinatorial distance at most $n$ from the central one may miss some elements of $\Dist_N$: for $N=9$, the $46$ distances obtained with $n=4$ range up to $7.3474$ but do not include the elements $7.2065$ and $7.2517$ of $\Dist_9$, which require larger $n$. A list for $N=9$ with no omissions up to $d_{14}$ was obtained there by a graphical procedure. We use instead the following simple observation, which gives a complete list up to any prescribed bound.

\begin{lemma}\label{lem:nogaps}
Let $D>0$, and let $S$ be the smallest set of centres of the tessellation with the following two properties: $O\in S$, and if $s\in S$ and $\rho(O,s)\le D+\rc$, then the centres of all the polygons sharing an edge with the polygon of $s$ belong to $S$. Then $S$ is finite and contains every centre at distance at most $D$ from $O$. Consequently $\{\rho(O,s)\mid s\in S,\ 0<\rho(O,s)\le D\}=\Dist_N\cap(0,D]$.
\end{lemma}
In practice $S$ is computed by exploring the tessellation from the central polygon outwards, adding at each step the neighbours of the polygons already found whose centres are at distance at most $D+\rc$ from $O$, until no new centre appears.
\begin{proof}
The set $S$ is finite because there are finitely many centres at distance at most $D+\rc$ from $O$, and each polygon has $N$ neighbours. Let $P$ be a centre with $\rho(O,P)\le D$. The geodesic segment $[O,P]$ crosses a finite sequence of polygons, from the polygon of $O$ to the polygon of $P$, consecutive ones sharing an edge or a vertex; since the three polygons meeting at a vertex share an edge pairwise, we may insert polygons so that consecutive ones always share an edge. Every polygon of the sequence contains a point at distance at most $D$ from $O$ (a point of the segment, or a vertex lying on it), so its centre is at distance at most $D+\rc$ from $O$. By induction along the sequence, all its centres belong to $S$; in particular $P\in S$.
\end{proof}

\subsection{The circumscribed circle}
The next statement reduces the search of all the packings of $X$ to a compact region of the fundamental domain. The use of the circumscribed disc for this purpose, and the fact that its radius $\rc$ cannot be lowered come from Chapter~7 of \cite{Reyes_thesis}.

\begin{lemma}\label{lem:rc}
Every closed disc of radius $\rc$ in $\D$ contains a centre of the tessellation $\{N,3\}$. In particular, every extremal packing of $X=\D/K$ has a disc whose centre lifts to the closed disc of radius $\rc$ about $O$.
\end{lemma}
\begin{proof}
The polygons of the tessellation are the Dirichlet--Voronoi cells of the set of centres, and every point of a cell is at distance at most $\rc$ from its centre. The second assertion follows by applying the first one to the lifted packing, whose Dirichlet--Voronoi cells form another regular tessellation $\{N,3\}$ with the same $\rc$.
\end{proof}

The value of $\rc$ is sharp: the open disc of radius $\rc$ about a vertex contains no centre.

\subsection{The number field}
As shown in Theorem~5.1 of \cite{Reyes_thesis}, the entries of the matrices of $a,b,c$ in $\PGL_2(\R)$ (acting on the upper half plane, orientation reversing elements acting on $\bar z$) can be chosen in the field $L_N=\Q(u,\cos\frac{\pi}{N})$ with $u=e^{R}$, namely
$$\bar a=\begin{pmatrix}1&0\\0&-1\end{pmatrix},\qquad \bar b=\begin{pmatrix}0&u^{-1}\\u&0\end{pmatrix},\qquad \bar c=\begin{pmatrix}-\cos\frac{\pi}{N}&\sin\frac{\pi}{N}\\ \sin\frac{\pi}{N}&\cos\frac{\pi}{N}\end{pmatrix},\qquad \sin\frac{\pi}{N}=\frac{u}{u^2+1}.$$
For $N\in\{7,8,9,12\}$ the cosine already lies in $\Q(u)$ and $L_N=\Q(u)$ has degree $12,8,12,8$; for $N=10$ and $N=18$ a quadratic extension is needed and $L_N$ has degree $8$ and $12$. We note that $L_9=L_{18}$ (the fields coincide, with $u_{18}\in\Q(u_9)$), a fact related to the inclusion of the centres of the $18$-gon tessellation among the centres of the $9$-gon one ($2R_{18}=d_2$ for $N=9$). All exact computations below take place in $L_N$ and were done with PARI/GP \cite{PARI2}.

\section{An algorithm for the determination of all the extremal packings}\label{sec:algorithm}

The input is the permutation action of $a,b,c$ on $\Dpm/K$. The output is the list of all the extremal packings of $X=\D/K$, each of them confirmed by an exact matrix $\alpha\in\PGL_2(L_N)$ with $\alpha^{-1}K\alpha<\Dpm$. The procedure has six steps.

\subsection*{Step 1: fundamental domain and side pairings}
The polygons of the tessellation modulo $K$ are the orbits of the dihedral group $\langle a,c\rangle$ on $\Dpm/K$; there are $k$ of them. We build a connected fundamental domain $F$ polygon by polygon: we start with the central polygon $P_0$, centred at $O$ and containing $T_0$, and, as long as $F$ has fewer than $k$ polygons, we add a polygon which is not $K$-equivalent to any polygon of $F$ and shares an edge with one of them (the polygon beyond edge $j$ of a polygon $g_P(P_0)$ of $F$ is $g_P(ca)^jab\,(P_0)$). In this way $F$ consists of $k$ polygons, each one glued to a previous one along an edge, and we obtain an explicit word $g_P\in\Dpm$ for each polygon $P$ of $F$ and for each triangle of $F$. Each edge of $F$ other than these gluing edges gives a side pairing $s\in K$, namely
$$s=g_P\,(ca)^j ab\,a^{\varepsilon}(ca)^{-m}\,g_Q^{-1},$$
where the triangle $g_P(ca)^jab\,(T_0)$ beyond edge $j$ of $P$ is $K$-equivalent to the triangle $g_Q(ca)^m a^{\varepsilon}(T_0)$ of $F$; the pairing is orientation reversing if and only if $\varepsilon=1$. The side pairings generate $K$ (Poincar\'e's theorem), and the fact that they fix the coset $K$ is checked on the permutation representation. The same element may pair two different pairs of edges, so the number of side pairings may be smaller than $(kN-2k+2)/2$ (for instance, nine instead of ten for the three octagons of the example with $N=8$ in Section~7.1 of \cite{Reyes_thesis}).

\subsection*{Step 2: admissible distances}
Using Lemma~\ref{lem:nogaps} we compute $\Dist_N\cap(0,D]$ for
$$D=2(\rc+\mu)+\max\{\rho(O,s_1O),\rho(O,s_2O),\,2R+\tfrac12\}+\tfrac1{10},$$
where $s_1,s_2$ are the two side pairings that move $O$ the least among those with different axes and $\mu=1/10$ is a safety margin discussed below; the last term $\frac1{10}$ is a further margin, independent of $\mu$. The first two terms in the maximum account for the displacements needed in Steps 3 and~4; the third one ensures that the list of centres computed here also contains those used in Step~6: since $\alpha^{-1}(O)$ is at distance at most $\rc+\mu$ from $O$, the centres at distance at most $\rc+2R+\frac{3}{10}$ from $\alpha^{-1}(O)$ are at distance less than $D$ from $O$. Repeated centres are detected \emph{exactly}: the coordinates of $g(O)$ are elements of $L_N(i)$, and for each distance the value $\tanh^2\frac12\rho(O,g(O))\in L_N$ is kept. For $N=9$ the resulting list coincides with the list $d_0,\dots,d_{14}$ obtained in Chapter~7 of \cite{Reyes_thesis}, and Lemma~\ref{lem:nogaps} guarantees that no element of $\Dist_9$ up to $d_{14}$ is missing from it.

\subsection*{Step 3: candidates}
Let $s$ be a hyperbolic isometry or a glide reflection, with translation length $T$, and let $d\in\Dist_N$. Following \cite{Girondo_Munoz}, we call
$$b_s(d)=\{z\in\D : \rho(z,sz)=d\}$$
an \emph{equi-displaced cycle}, or \emph{ED-cycle}. It is the union of the two hypercycles at distance $\eta$ from the axis of $s$, where
$$\cosh\eta=\frac{\sinh(d/2)}{\sinh(T/2)}\ \text{ if $s$ is hyperbolic},\qquad \cosh\eta=\frac{\cosh(d/2)}{\cosh(T/2)}\ \text{ if $s$ is a glide reflection};$$
it reduces to the axis when the right-hand side equals $1$, and it is empty when it is smaller. In the disc model each hypercycle is a circular arc through the end points of the axis (see \cite[Theorem 7.35.1]{Beardon_1983} for hyperbolic elements and \cite[Theorems~4.4, 4.5 and Corollary~4.6.1]{Reyes_thesis} for glide reflections, where these sets are called \emph{bananas}).

By Lemma~\ref{lem:rc} the centre of some disc of any extremal packing lies in the closed disc $\overline{\D}(O,\rc)$, where its displacement by $s_i$ is at most $2\rc+\rho(O,s_iO)$. Hence every such centre is one of the finitely many intersection points of the ED-cycles $b_{s_1}(d)\cap b_{s_2}(d')$ with $d,d'\in\Dist_N$ bounded as above (the finiteness of the number of extremal packings of a surface, which follows from the same argument, is \cite[Corollary~4.6.2]{Reyes_thesis}). We compute all of them, keep those in $\overline{\D}(O,\rc+\mu)$ --- the margin $\mu$ takes care of tangential intersections on the boundary circle, which do occur --- and discard those which are moved by some side pairing a distance that is \emph{not} in $\Dist_N$ (how these decisions are made rigorous is explained after Proposition~\ref{prop:cert}). The survivors are the \emph{candidates}; typically there are between two and ten of them, the centre $O$ itself among them.

\subsection*{Step 4: numerical search of $\alpha$}
Let $z$ be a candidate, and suppose that it is the centre of a disc of a packing with tessellation $\alpha\Dpm(T_0)$. Composing $\alpha$ with an element of $\Dpm$ we may assume that $\alpha(O)=z$ and that $\alpha$ preserves the orientation. Let $\beta$ be the translation along the geodesic through $O$ and $z$ that maps $O$ to $z$ (the identity if $z=O$). Then $\beta^{-1}\alpha$ fixes $O$, so it is a rotation $R_\theta$ about $O$, and $\alpha=\beta R_\theta$. To find $\theta$, note that $s_1z$, where $s_1\in K$ is the side pairing chosen in Step~2, is also a centre of the tessellation: $s_1z=\alpha(g(O))$ for some $g\in\Dpm$ with $\rho(O,g(O))=\rho(z,s_1z)$. Hence $\beta^{-1}(s_1z)=R_\theta(g(O))$, and since there are finitely many centres $g(O)$ at that distance from $O$, this gives finitely many angles $\theta$, modulo $2\pi/N$. For each of them we test whether $\alpha=\beta R_\theta$ conjugates every side pairing into $\Dpm$. The test is a \emph{reduction algorithm}. Fix an interior point $q_0$ of $T_0$ and let $h$ be an isometry. As long as the point $p=h(q_0)$ is not in $T_0$, we replace $p$ by its image under the reflection in a side of $T_0$ whose line separates $p$ from $T_0$. Each such reflection decreases by one the number of lines of the tessellation $\Dpm(T_0)$ separating $p$ from $T_0$, as in any discrete group generated by reflections, so after finitely many steps $p\in T_0$, and the reflections used, read in order, form a word $w$ in $a,b,c$ with $w^{-1}h(q_0)\in T_0$. If $h\in\Dpm$, then $w^{-1}h$ is an element of $\Dpm$ mapping $T_0$ to itself, hence the identity; so $h\in\Dpm$ if and only if $w^{-1}h=1$, and in that case $h=w$. Applied to $h=\alpha^{-1}s\alpha$ for every side pairing $s$, this produces, when $\alpha$ is correct, words $w_s$ with $\alpha^{-1}s\alpha=w_s$; with the approximate $\alpha$ these equalities hold only approximately, and they are made exact in Step 5.

\subsection*{Step 5: exact verification}
The numerical $\alpha$ and the words $w_s$ are turned into a proof as follows. Let $\bar s$ and $\bar w_s$ be the exact matrices in $\PGL_2(L_N)$ of the side pairings and of the words. The conditions $\bar s\,\alpha=\pm\alpha\,\bar w_s$, one for each side pairing (the signs being read off from the numerical solution), form a homogeneous linear system over $L_N$ in the four entries of $\alpha$; we compute its kernel with PARI.

\begin{proposition}\label{prop:cert}
If the kernel is one-dimensional, spanned by a matrix $\alpha$ with $\det\alpha\neq0$, then $\alpha^{-1}K\alpha<\Dpm$ and the points $\alpha g(O)$, $g\in\Dpm$, are the lifted centres of an extremal $k$-packing of $X$. Conversely, every extremal packing of $X$ has a disc centred at one of the candidates of Step~3.
\end{proposition}
\begin{proof}
The equations say that $\alpha^{-1}\bar s\alpha=\pm\bar w_s\in\Dpm$ for a generating set of $K$, hence $\alpha^{-1}K\alpha<\Dpm$, i.e.\ $K<\alpha\Dpm\alpha^{-1}$ preserves the regular tessellation $\alpha\Dpm(T_0)$ and the B\"or\"oczky packing of its inscribed discs, which therefore descends to a packing of $X$ by $k$ discs of radius $R$. The converse is Lemma~\ref{lem:rc} together with the exhaustive enumeration of Step 3. Thus, once every candidate has either been confirmed in this way or shown not to be a centre (see below), the list of packings is complete.
\end{proof}

Note that this verification involves neither the coordinates of $z$ nor any square root: the whole verification takes place in $L_N$. In Chapter~5 of \cite{Reyes_thesis} and in \cite{Girondo_Reyes_2023} the algebraic verification was carried out by constructing explicitly an automorphism of the surface (a rotation of order two) from the exact coordinates of the candidate, which required a further quadratic extension of $L_N$ for each surface; the linear formulation of Proposition~\ref{prop:cert} avoids both. The centres of the packings found are the points $\alpha g(O)$, which have coordinates in $L_N(i)$, so they are known exactly; decimal values are given below only for orientation. For the surface $K_{(9,-,2)}$ studied in detail in Section~7.1 of \cite{Reyes_thesis} the hidden centre found by the algorithm, approximately $0.447558+0.016583i$, is the one obtained there.

\subsection*{Exact and floating point steps}
Steps 1, 2 and 5 are exact: they use only permutations, words in $a,b,c$, and arithmetic in $L_N$ and $L_N(i)$. Steps 3 and 4 use floating point arithmetic (with $40$ significant digits), and their results enter the proof in two different ways.

A \emph{positive} outcome of Step 4, a candidate $z$ with a conjugator $\alpha$, is never used as such: it is only the input of Step 5, which either proves the existence of the packing (Proposition~\ref{prop:cert}) or fails. For all the surfaces of this paper Step 5 never failed: every conjugator found in Step 4 gave a one-dimensional kernel spanned by an invertible matrix.

A \emph{negative} outcome, a point discarded in Step 3 or a candidate for which Step 4 finds no conjugator, is justified by the basic property of admissible distances (Section~\ref{sec:observations}; \cite[Lemma~4.1]{Reyes_thesis}): if $z$ is the centre of a disc of an extremal packing, then $\rho(z,\kappa z)\in\Dist_N$ for every $\kappa\in K$. A point is discarded in Step 3 because $\rho(z,sz)\notin\Dist_N$ for some side pairing $s$; for every candidate without a conjugator we exhibit in the same way an element $\kappa=s s'^{\pm1}\neq1$, product of two side pairings, with $\rho(z,\kappa z)\notin\Dist_N$, so that it is not a centre either (for $10$ of these $70$ candidates, on two surfaces with $N=7$, the displacement $\rho(z,\kappa z)$ exceeds the bound $D$ of Step~2, and the decision uses the complete list of admissible distances up to a slightly larger bound, given by Lemma~\ref{lem:nogaps}). These are the only decisions in which floating point numbers intervene, and each of them is taken with a wide margin. Over all the surfaces of this paper and their double covers, an equality $\rho(z,sz)=d$ with $d\in\Dist_N$ is reproduced by the floating point computation with an error of at most $2\cdot10^{-19}$ (the largest errors occur at nearly tangential intersections, where the intersection point is obtained from the square root of a number close to $0$), whereas for every discarded point $z$ there is a side pairing $s$ with $|\rho(z,sz)-d|\geq 6\cdot10^{-6}$ for every $d\in\Dist_N$, and for every candidate without conjugator there is a product $\kappa\neq1$ of two side pairings with $|\rho(z,\kappa z)-d|\geq 2.7\cdot10^{-3}$ for every $d\in\Dist_N$. The threshold used by the program to decide whether a displacement belongs to $\Dist_N$ is $10^{-8}$, so the rounding errors are more than ten orders of magnitude too small to change any of these decisions. These bounds are recomputed by the program \texttt{margins.py} (Section~\ref{sec:data}).

\subsection*{Step 6: counting packings}
Two verified pairs $(z,\alpha)$, $(z',\alpha')$ give the same packing of $X$ if and only if the lifted centre sets $K\alpha\Dpm(O)=\alpha\Dpm(O)$ and $\alpha'\Dpm(O)$ coincide, and since a regular tessellation is determined by one centre together with its $N$ neighbours this can be decided by comparing the centres lying at distance at most $\rc+2R+\frac{3}{10}$ from $O$. The number of these classes is the number of extremal $k$-packings of $X$. It may happen that two different packings share a disc (two of the five $4$-packings of the surface $X_6$ of Theorem~\ref{thm:X6} share the disc centred at $O$); this is why packings, and not centres, are counted.

\subsection{Automorphisms, real structures and hidden isomorphisms}\label{sec:hidden}
The algorithm also gives the full isometry group of $X$. Indeed, an isometry $\gamma$ of $X$ lifts to an element of the normalizer of $K$ in $\Isom(\D)$ and maps the packing $\mathcal P_0$ with centres $\Dpm(O)$ to one of the packings $\mathcal P_j$ with centres $\alpha_j\Dpm(O)$; since $\Dpm$ is the full symmetry group of the tessellation, $\gamma\in\alpha_j\Dpm$. Writing $K_j=\alpha_j^{-1}K\alpha_j<\Dpm$, which is explicitly known through the words $w_s$, we get:
\begin{itemize}[leftmargin=2em]
\item the isometries of $X$ fixing $\mathcal P_0$ are the \emph{induced} ones, $N_{\Dpm}(K)/K$;
\item there are isometries of $X$ mapping $\mathcal P_0$ to $\mathcal P_j$ if and only if $K_j$ is conjugate to $K$ in $\Dpm$ (if and only if the permutation representations of $\Dpm$ on $\Dpm/K$ and on $\Dpm/K_j$ are isomorphic), and then there are exactly $|N_{\Dpm}(K)/K|$ of them, $\alpha_j\delta$ with $\delta K\delta^{-1}=K_j$;
\item for orientable $X$, such an isometry is anticonformal if and only if $K_j$ is $\Dl$-conjugate to $\sigma K\sigma^{-1}$ ($\sigma\in\Dpm\smallsetminus\Dl$, as in Section~\ref{sec:classes}) rather than to $K$;
\item if $K_j$ is $\Dpm$-conjugate to a representative $K'$ of a \emph{different} class, then $\D/K$ and $\D/K'$ are isometric, although $K$ and $K'$ are not conjugate in the triangle group.
\end{itemize}
The second item has a consequence that explains when two extremal packings of the same surface are not related by an isometry.

\begin{lemma}\label{lem:orbits}
Let $X=\D/K$ with $K<\Dpm$. The orbits of $\Isom(X)$ on the set of extremal packings of $X$ are in bijection with the $\Dpm$-conjugacy classes of the subgroups $K'<\Dpm$ such that $\D/K'$ is isometric to $X$, the orbit of $\mathcal P_j$ corresponding to the class of $K_j$. In particular, $X$ carries two extremal packings that are not related by an isometry if and only if $X$ is uniformized by two subgroups of $\Dpm$ that are not conjugate in $\Dpm$.
\end{lemma}

\begin{proof}
By the second item, $\mathcal P_i$ and $\mathcal P_j$ are in the same orbit if and only if $K_i$ and $K_j$ are conjugate in $\Dpm$ (apply it with $\mathcal P_i$ in the role of $\mathcal P_0$, i.e.\ to $X=\D/K_i$), so the map is well defined and injective; each $K_j=\alpha_j^{-1}K\alpha_j$ uniformizes a surface isometric to $X$. Conversely, if $K'<\Dpm$ and $\D/K'$ is isometric to $X$, there is $\varphi\in\Isom(\D)$ with $\varphi K\varphi^{-1}=K'$; then $K<\varphi^{-1}\Dpm\varphi$, so $\alpha=\varphi^{-1}$ defines an extremal packing of $X$ whose group $\alpha^{-1}K\alpha$ is $K'$.
\end{proof}

Comparing transitive permutation representations is done through a canonical form, so all these questions are decided exactly. Orders and fixed points of the resulting (possibly hidden) automorphisms are computed by reducing powers into $\Dpm$ and by searching the elliptic elements of $K\gamma$ whose fixed points lie in $F$.

\subsection{Implementation and validation}\label{sec:impl}
The whole procedure is implemented in Python (floating point arithmetic with \texttt{mpmath} \cite{mpmath}, exact arithmetic with PARI/GP through \texttt{cypari2}), with GAP \cite{GAP4} for the enumeration of the groups and for the coset actions. On a standard computer one surface takes between one and three seconds, and the complete run over the $373$ surfaces and the $50$ canonical double covers about five minutes. The results agree with all the previous classifications. For the surface $K_{(9,-,2)}$ worked out in detail in Section~7.1 of \cite{Reyes_thesis}, the permutations $a,b,c$ printed there define exactly our class \texttt{nec\_N9\#2}, which confirms that the conventions for the coset action agree, and the program reproduces the list of admissible distances, a fundamental domain with eight side pairings as there, the hidden centre and the two packings of the surface and of its double cover. The number of packings of each of the $50$ non-orientable surfaces and of their double covers agrees with \cite[Theorem~7.2]{Reyes_thesis}. The algorithm is not restricted to the families of this paper: run on the $144$ torsion-free NEC subgroups of index $36$ of $\Dpm(2,3,18)$, that is, on the non-orientable extremal surfaces of genus $4$ classified in \cite{Nakamura_2009}, it recovers the $22$ surfaces with two extremal discs, the five possible positions of the second centre (at the points $\tanh(R/2)$, $1/\sinh R$, $1/(2\tanh R\cos\beta)$, $1/(2\tanh R\cos 2\beta)$ and $2\sin 2\beta/\tanh R$ of the disc, $\beta=\pi/18$) and the orders of the automorphism groups of Tables~1 and~2 of that paper, surface by surface. All the data (groups, side pairings, exact matrices $\alpha$, words, automorphisms) are available in the files described in Section~\ref{sec:data}.

\section{A worked example: the surface $X_6$}\label{sec:example}

\begin{figure}[h]
\centering
\includegraphics[width=0.62\textwidth]{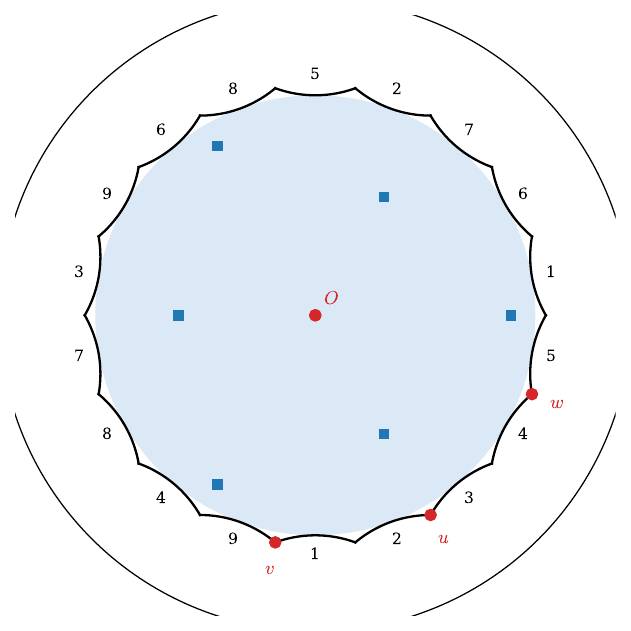}
\caption{The fundamental $18$-gon of $X_6$ (edges with the same label are identified), the extremal disc inscribed in it, the centres of the four extremal discs (red: the centre $O$ and the three vertices $u,v,w$) and the six Weierstrass points (blue squares).}\label{fig:X6domain}
\end{figure}

We illustrate the six steps of the algorithm on the surface that will turn out to be the most interesting one, the surface uniformized by a torsion-free subgroup $K_{18}$ of index $18$ in $\Dl(2,3,18)$ (\texttt{fuchs\_N18\#9}), that is, a $1$-extremal Riemann surface of genus $2$ with an $18$-gon as fundamental domain (Figure~\ref{fig:X6domain}). It is the surface $X_6$ of \cite{Girondo_Gonzalez-Diez_2002_2}, where its four extremal discs, its automorphism group and its Weierstrass points were determined by a geometric analysis of the $18$-gon; the example shows how the algorithm recovers these facts mechanically. Its group $K$ is generated by the four elements
$$(bcaca)^2,\qquad bcabcbabcabcba,\qquad cb(ac)^2babcba,\qquad cbabcbacbaca$$
(the generators returned by GAP, freely reduced using $a^2=b^2=c^2=1$),
and its permutation representation on the $36$ cosets of $\Dpm(2,3,18)$ is the input of the algorithm.

\emph{Step 1.} There is a single polygon, so no gluing is needed and the nine side pairings are the elements
$$s_j=(ca)^{j}ab\,(ca)^{-m_j}\in K\qquad\text{pairing edge $j$ with edge $m_j$,}$$
for the nine pairs $\{j,m_j\}$ of Figure~\ref{fig:X6domain}; in terms of the edge labels $0,\dots,17$ (edge $0$ at the bottom, counterclockwise) they are $\{0,5\}$, $\{1,8\}$, $\{2,13\}$, $\{3,16\}$, $\{4,9\}$, $\{6,11\}$, $\{7,14\}$, $\{10,15\}$, $\{12,17\}$. All of them are hyperbolic translations moving the centre $O$ a distance $d_0=2R\approx3.438214$.

\emph{Step 2.} The list of admissible distances for $N=18$ is computed up to $D\approx7.95$ and contains $16$ values, the first of which are those of Table~\ref{ta:dist}.

\begin{figure}[h]
\centering
\includegraphics[width=0.6\textwidth]{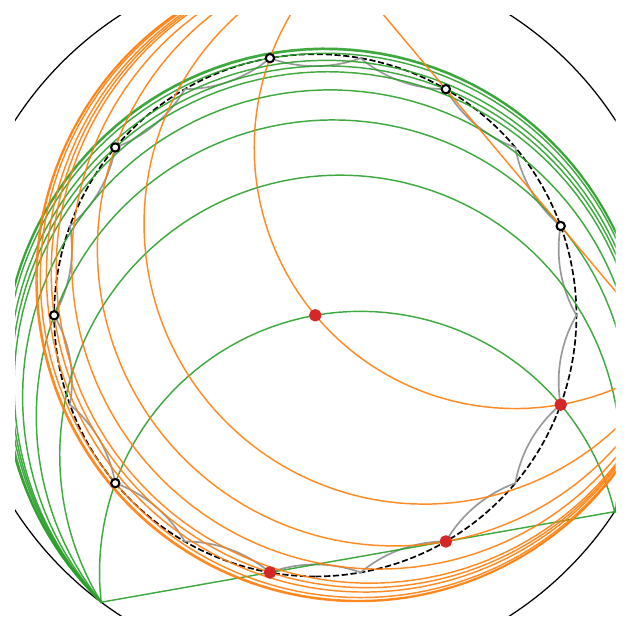}
\caption{Step 3 for $X_6$: the hypercycles $\{\rho(z,s_1z)=d\}$ (green) and $\{\rho(z,s_2z)=d\}$ (orange), $d\in\Dist_{18}$, meeting the search disc (dashed). The candidates are the intersection points marked with circles; the centres confirmed in Step 5 are the red ones.}\label{fig:edcycles}
\end{figure}

\emph{Step 3.} Two side pairings with different axes are chosen, say $s_1$ and $s_2$ (pairing the edges $\{0,5\}$ and $\{1,8\}$). The ED-cycles $b_{s_i}(d)$ for $d\in\Dist_{18}$, $d\le 2(\rc+\mu)+d_0$, are drawn in Figure~\ref{fig:edcycles}; they meet at $61$ points inside the search disc of radius $\rc+\mu$ about $O$, and the displacement test with respect to all nine side pairings leaves $10$ candidates: the centre $O$ and nine vertices of the polygon (the vertices are at distance $\rc$ from $O$, so this is a situation where the margin $\mu$ matters: the intersections of the hypercycles at the vertices are tangential).

\emph{Steps 4--6.} For the candidate $O$ the only admissible rotation is the identity (the original disc); for three of the nine vertices a conjugator $\alpha=\beta R_\theta$ is found, all its $\alpha^{-1}s_j\alpha$ reduce to words of length between $22$ and $50$ in $a,b,c$, and the linear system of Step 5 has a one-dimensional kernel, which confirms in exact arithmetic the three hidden extremal discs centred at the vertices $u,v,w$ of Figure~\ref{fig:X6domain}; the other six vertices are rejected (some product of two side pairings moves them a non-admissible distance). The four lifted centre sets are pairwise different, so $X_6$ has exactly four extremal discs, in agreement with \cite{Girondo_Gonzalez-Diez_2002_2}.

\emph{Automorphisms.} The group induced by $\Dpm(2,3,18)$ has order $6$: the rotation of order $3$ about $O$ and three reflections. For each of the three hidden discs the group $K_j=\alpha_j^{-1}K\alpha_j$ is $\Dpm$-conjugate to $K$, so each hidden disc contributes six more isometries and $|\Isom(X_6)|=24$, $|\Aut^+(X_6)|=12$, as in \cite{Girondo_Gonzalez-Diez_2002_2}; the explicit hidden isometries are the products $\alpha_j\delta$ of Section~\ref{sec:hidden}. Among the hidden conformal involutions exactly one has six fixed points: it is the hyperelliptic involution, and its fixed points, computed as fixed points of the elliptic elements $\kappa\alpha_j\delta$, $\kappa\in K$, are the Weierstrass points shown in Figure~\ref{fig:X6domain}, which are those found in Proposition~1 of \cite{Girondo_Gonzalez-Diez_2002_2}.

\begin{figure}[h]
\centering
\includegraphics[width=0.95\textwidth]{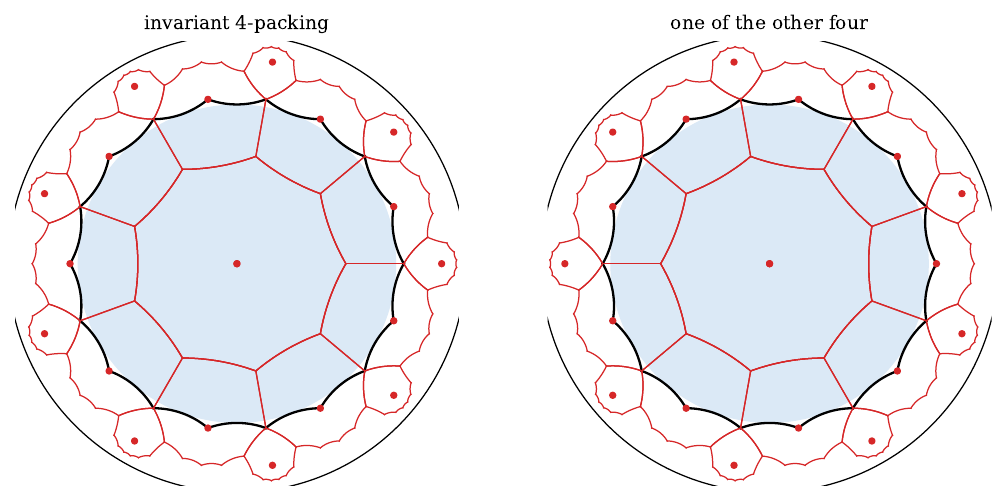}
\caption{Two of the five extremal $4$-packings of $X_6$ drawn in the $18$-gon: the $9$-gons of the packing invariant under the whole isometry group (left) and of one of the four packings forming an orbit (right). Both have a disc centred at $O$.}\label{fig:twotess}
\end{figure}

\emph{The same surface as a $4$-extremal surface.} Since $L_9=L_{18}$, the exact verification of Step 5 can be carried out for the group $K$ with respect to $\Dl(2,3,9)$: the candidates are now the points of the disc of radius $\rc(9)\approx1.035886$ about $O$ whose displacements under the $s_j$ lie in $\Dist_9$, and the search of $\alpha$ is made with the centres of the $\{9,3\}$ tessellation. It finds five conjugates $\alpha\Dl(2,3,9)\alpha^{-1}$ containing $K$, i.e.\ five extremal $4$-packings; two of them have a disc centred at $O$ and differ by a rotation of angle $\pi/9$ about $O$ (Figure~\ref{fig:twotess}), the other three have their nearest centre at distance exactly $\rc(9)$ from $O$. The groups $\alpha^{-1}K\alpha<\Dl(2,3,9)$ fall into two $\Dl(2,3,9)$-conjugacy classes, those of the groups $K_9'$ (four of the packings) and $K_9$ (one packing) of Theorem~\ref{thm:X6}, which is how the identification of Theorem~\ref{thm:X6} is obtained.

\section{Results}\label{sec:results}

\subsection{Non-orientable surfaces of genus 3}

\begin{theorem}\label{thm:nonor}
The number of extremal $k$-packings of the non-orientable $k$-extremal surfaces of genus $3$ is given by the following table, where $m(n)$ means that $n$ surfaces admit exactly $m$ packings.
\begin{center}
\begin{tabular}{cccl}
\toprule
$k$ & $N$ & number of surfaces & extremal $k$-packings\\
\midrule
$6$ & $7$ & $12$ & $1\,(3)$, $2\,(9)$\\
$3$ & $8$ & $23$ & $1\,(8)$, $2\,(15)$\\
$2$ & $9$ & $4$ & $1\,(1)$, $2\,(3)$\\
$1$ & $12$ & $11$ & $1\,(1)$, $2\,(10)$\\
\bottomrule
\end{tabular}
\end{center}
No two of these $50$ groups are conjugate in $\Isom(\D)$, so the $50$ surfaces are pairwise non-isometric. In every surface with two packings, the two packings are exchanged by an isometry of the surface.
\end{theorem}

The groups of isometries induced by $\Dpm$ have orders $1$, $2$, $4$, $6$ or $8$. The surfaces with an induced group of order $4$ or $8$ (they only occur for $N=8$) have a single packing, all those with trivial induced group or with induced group of order $6$ have two packings, and for induced group of order $2$ both possibilities occur. (The $11$ surfaces with $k=1$ had been studied in \cite{Girondo_Nakamura_2007}; the whole table was obtained in \cite{Reyes_thesis} by a surface-by-surface analysis.)

\subsection{Double covers}

\begin{theorem}\label{thm:cover}
Let $X$ be a non-orientable $k$-extremal surface of genus $3$ and $X^+$ its canonical double cover, a $2k$-extremal Riemann surface of genus $2$. Then $X^+$ has exactly as many extremal $2k$-packings as $X$ has extremal $k$-packings. The group of isometries of $X^+$ induced by the triangle group contains the lift of the one induced on $X$, which has twice its order, and coincides with it except for four surfaces (\texttt{nec\_N12\#3}, \texttt{nec\_N8\#4}, \texttt{nec\_N8\#12} and \texttt{nec\_N7\#4}), for which the group induced on $X$ has order $2$ and the one induced on $X^+$ has order $12$. The $50$ double covers are pairwise non-isomorphic; they are the $11$, $4$, $23$ and $12$ classes of Fuchsian subgroups of $\Dl(2,3,N)$, $N=12,9,8,7$, listed in Appendix~\ref{app:covers}.
\end{theorem}

The double cover $K^+=K\cap\Dl$ is handled by the algorithm through its coset action, obtained from that of $K$ by doubling the points with a parity label. In general the double cover has at least as many extremal packings as the non-orientable surface (\cite[Theorem~6.2]{Reyes_thesis}; see also \cite{Girondo_Reyes_2023}); the equality in genus $3$ was obtained in \cite[Theorem~7.2]{Reyes_thesis} by the surface-by-surface analysis, and is confirmed here by an independent method; this equality is special to genus~$3$: it fails for five of the non-orientable extremal surfaces of genus~$4$ (Remark~\ref{rem:genus4}). The statement about the induced isometry groups is new.

Among the Fuchsian classes with an anticonformal automorphism induced by $\Dpm$ (column 6 of Table~\ref{ta:all_info}), those which are \emph{not} canonical double covers --- $4$, $7$, $14$ and $13$ classes for $N=12,9,8,7$ --- are exactly those whose induced anticonformal automorphisms are either reflections (involutions with fixed points) or have order $4$, $8$ or $12$ without fixed points; none of them admits a fixed point free anticonformal involution induced by the triangle group, as it must be. For instance the $6$-extremal Riemann surface of genus $2$ uniformized by the normal subgroup of index $48$ of $\Dl(2,3,8)$, i.e.\ Bolza's surface (Remark~\ref{rem:bolza} and Figure~\ref{fig:classical}), has $96$ induced isometries, a single extremal $6$-packing, and $18$ induced reflections.

\subsection{Orientable surfaces of genus 2}

The line $k=1$ of the table in the following theorem recovers the classification of the $1$-extremal surfaces of genus $2$ obtained in \cite{Girondo_Gonzalez-Diez_2002_2, Nakamura_2002} (see Remark~\ref{rem:ggd}); for $k=2,3,4,6,12$ the $50$ surfaces that are canonical double covers of non-orientable ones had been classified in \cite{Reyes_thesis} (see Theorem~\ref{thm:cover}), and some further examples with two packings had been given in \cite{Girondo_Reyes_2023} and in Section~6.1 of \cite{Reyes_thesis}; the rest of the classification is new.

\begin{theorem}\label{thm:or}
The number of extremal $k$-packings of the $k$-extremal Riemann surfaces of genus $2$ is given by the following table, where surfaces are counted up to isometry (conjugacy classes in $\Dpm$, identified as in the last column).
\begin{center}
\begin{tabular}{ccccl}
\toprule
$k$ & $N$ & classes in $\Dpm$ & surfaces & extremal $k$-packings\\
\midrule
$1$ & $18$ & $8$ & $8$ & $1\,(1)$, $2\,(6)$, $4\,(1)$\\
$2$ & $12$ & $20$ & $20$ & $1\,(5)$, $2\,(15)$\\
$3$ & $10$ & $15$ & $15$ & $1\,(4)$, $2\,(11)$\\
$4$ & $9$ & $24$ & $23$ & $1\,(8)$, $2\,(13)$, $4\,(1)$, $5\,(1)$\\
$6$ & $8$ & $57$ & $57$ & $1\,(33)$, $2\,(24)$\\
$12$ & $7$ & $90$ & $90$ & $1\,(60)$, $2\,(30)$\\
\bottomrule
\end{tabular}
\end{center}
Moreover:
\begin{enumerate}[label=(\roman*),leftmargin=2.5em]
\item Exactly one pair of $\Dpm$-conjugacy classes gives rise to the same surface: two of the $24$ classes of $4$-extremal surfaces ($N=9$), both with five packings, uniformize the surface $X_6$ of Theorem~\ref{thm:X6}. In the table $X_6$ has been counted once (in the line $N=9$); in the line $N=18$ it is the surface with four extremal discs.
\item A $k$-extremal Riemann surface of genus $2$ admits an anticonformal automorphism if and only if it admits one induced by the triangle group. Consequently the numbers of $k$-extremal Riemann surfaces up to conformal equivalence are $9$, $25$, $20$, $36$, $77$ and $155$ for $k=1,2,3,4,6,12$, the mirror images of the surfaces counted by $r$ in Section~\ref{sec:classes} being pairwise non-isomorphic.
\item Only two surfaces have more than two packings: $X_6$ and the $4$-extremal surface $X_3$ of Remark~\ref{rem:21}, with four packings. Two different extremal packings of the same surface share a disc only in the three surfaces of (v): in $X_2$ the two $4$-packings share one disc, in $X_3$ the four $4$-packings form two pairs and the two packings of each pair share one disc, and in $X_6$ the $4$-packing invariant under the isometry group shares exactly one disc with each of the other four.
\item The full conformal automorphism groups of the $23$ four-extremal surfaces are $\Z_2$ (six surfaces), $\Z_2\times\Z_2$ (sixteen) and $D_6$ (only $X_6$).
\item Exactly three of the eight $1$-extremal surfaces of genus $2$ are also $4$-extremal, namely the surfaces $X_2$, $X_3$ and $X_6$ of \cite{Girondo_Gonzalez-Diez_2002_2}: $X_2$ has one extremal disc and two extremal $4$-packings, $X_3$ has two extremal discs and four extremal $4$-packings (Remark~\ref{rem:21}), and $X_6$ has four discs and five $4$-packings. By Theorem 5 of \cite{Girondo_2018} these are all the surfaces of genus $2$ which are $k$-extremal for two different values of $k$.
\item In every surface other than $X_6$ the isometry group acts transitively on the set of extremal packings. In $X_6$ the five $4$-packings form two orbits, of sizes $4$ and $1$: there is an extremal $4$-packing of $X_6$ which is not the image of any of the other four under an isometry of the surface (the four extremal discs of $X_6$, on the other hand, form a single orbit).
\end{enumerate}
\end{theorem}

The second and last items are obtained by testing, for every hidden packing $\mathcal P_j$ of every surface, whether $K_j$ is $\Dpm$-conjugate to $K$ and, if so, whether it is $\Dl$-conjugate to $K$ or to $\sigma K\sigma^{-1}$ (Section~\ref{sec:hidden}): $K_j$ is always $\Dpm$-conjugate to $K$ except for the hidden packings of $X_6$ relating its two orbits, and it is never $\Dl$-conjugate to $\sigma K\sigma^{-1}$. By Lemma~\ref{lem:orbits}, items (i) and (vi) are two faces of the same fact: $X_6$ is the only surface of our list with two orbits of extremal packings because it is the only one uniformized by two subgroups of the same triangle group that are not conjugate in it. The same test applied to the non-orientable surfaces gives the last assertion of Theorem~\ref{thm:nonor}. In the same test $K_j$ is frequently $\Dl$-conjugate to $K$ itself, producing conformal automorphisms that are not induced by the triangle group; for example the surface $X_3$ of Theorem~\ref{thm:or}(iii) has a trivial induced conformal group and full conformal automorphism group of order $4$ (see Remark~\ref{rem:21}).

\begin{figure}[h]
\centering
\includegraphics[width=0.95\textwidth]{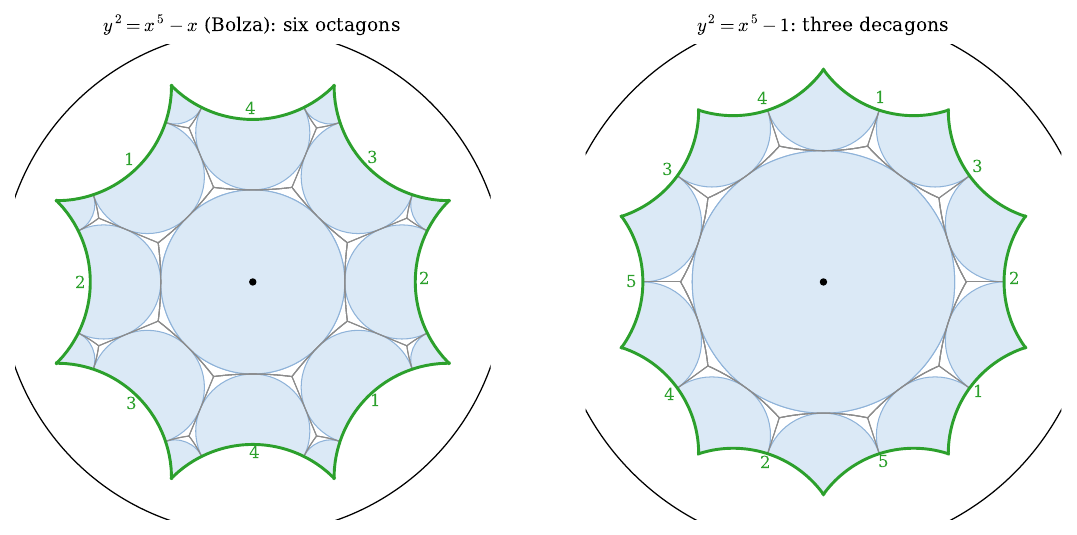}
\caption{The surfaces of Remark~\ref{rem:bolza}: Bolza's surface $y^2=x^5-x$ as a $6$-extremal surface ($N=8$, \texttt{fuchs\_N8\#77}) and the curve $y^2=x^5-1$ as a $3$-extremal surface ($N=10$, \texttt{fuchs\_N10\#18}). In green, the Dirichlet domain centred at the centre of a disc of the extremal packing (sides with the same number are identified); inside it, the tessellation by regular $N$-gons with angles $2\pi/3$ and the discs of the packing inscribed in them, cut along the sides of the domain. The six octagons and the three decagons of the packings can be reassembled from the pieces.}\label{fig:classical}
\end{figure}

\begin{remark}\label{rem:bolza}
Among the surfaces with a large induced group we find the classical curves: Bolza's surface $y^2=x^5-x$ as a $6$-extremal surface (induced group of order $96$, $N=8$) and the curve $y^2=x^5-1$, the only genus two curve with an automorphism of order $5$, as a $3$-extremal surface (induced group of order $20$, $N=10$). Both have a single extremal packing. Their Dirichlet domains centred at a disc of the packing are the classical regular polygons: an octagon with opposite sides identified for Bolza's surface and a decagon for $y^2=x^5-1$ (Figure~\ref{fig:classical}).
\end{remark}

\subsection{The surface $X_6$}
It is shown in \cite{Girondo_Gonzalez-Diez_2002_2} that the surface $X_6$ has exactly four extremal discs, centred at the points $o,u,v,w$ of the fundamental $18$-gon, that $\Aut^+(X_6)\cong D_6$ acts transitively on them, that $X_6$ has an anticonformal involution $S$ induced by a symmetry of the $18$-gon (so that its isometry group has order $24$), that $X_6$ is a curve $y^2=(x^3-1)(x^3-\lambda^6)$ for some $\lambda$ with $\lambda^6\neq1$ (proof of Lemma~1), and its six Weierstrass points are located there as fixed points of an explicit lift of the hyperelliptic involution $J$, which maps $o$ to $w$ (Proposition~1). The surface was pointed out again in Section 4.5 of \cite{Girondo_2018} as an example of a surface which is both $1$-extremal and $4$-extremal and contains two extremal $4$-packings not related by an automorphism. The following theorem completes its description; the new information is the number and the orbit structure of the $4$-packings, the identification of the uniformizing groups, the value of $\lambda$, and the description of the Weierstrass points in terms of the tessellation of a $4$-packing.

\begin{theorem}\label{thm:X6}
The surface $X_6$ of \cite{Girondo_Gonzalez-Diez_2002_2}, which is $1$-extremal with exactly four extremal discs, is $4$-extremal with exactly five extremal $4$-packings. Up to conjugation in $\PSL_2(\R)$, its uniformizing group is the subgroup $K_{18}$ of index $18$ in $\Dl(2,3,18)$ of Section~\ref{sec:example}, and also each of two subgroups of index $36$ in $\Dl(2,3,9)$ which are not conjugate in $\Dpm(2,3,9)$: $K_9$, whose normalizer in $\Dpm(2,3,9)$ has index $3$, and $K_9'$, whose normalizer has index $12$ (\texttt{fuchs\_N9\#37} and \texttt{fuchs\_N9\#36}). Moreover:
\begin{enumerate}[label=(\roman*),leftmargin=2.5em]
\item ($\Aut^+(X_6)\cong D_6$ and $|\Isom(X_6)|=24$ by \cite{Girondo_Gonzalez-Diez_2002_2}.) Every isometry of $X_6$ preserves the decomposition of $X_6$ into the four $9$-gons of one of its extremal $4$-packings, which is therefore invariant under the whole isometry group (in the uniformization by $K_9$ this decomposition comes from the standard tessellation, and every isometry is induced by $\Dpm(2,3,9)$). The other four $4$-packings form a single orbit, as do the four extremal discs; each of these four packings shares exactly one disc with the invariant one. Thus the invariant $4$-packing is not related to the other four by any isometry: the two kinds of $4$-packings of $X_6$ are geometrically different. In terms of Lemma~\ref{lem:orbits}, the two orbits correspond to the two non-conjugate groups: the invariant packing is the standard one for $K_9$, and each packing of the orbit is the standard one for a conjugate of $K_9'$.
\item $X_6$ is the curve
$$y^2=x^6+2\sqrt{3}\,i\,x^3+1\ \cong\ y^2=(x^3-1)(x^3+(2+\sqrt3)^2)=(x^3-1)(x^3+7+4\sqrt3).$$
\item The hyperelliptic involution preserves the decomposition of $X_6$ given by the invariant $4$-packing (in the uniformization by $K_9$ it is induced by $ab$), and the six Weierstrass points are midpoints of edges of this decomposition. In the lifted tessellation, for a suitable $9$-gon $P$ (in the uniformization by $K_9$, the one centred at $O$), they are the midpoints of three alternate edges of $P$, at distance $R_9$ from its centre, and the midpoints of the three edges that leave $P$ at the vertices opposite to them, at distance $d_W$ from its centre, where
$$\cosh R_9=\frac{1}{2\sin(\pi/9)},\qquad \cosh d_W=\cosh 2R_9\cosh R_9-\sinh 2R_9\sinh R_9\cos\tfrac{2\pi}{9}$$
($R_9\approx0.927539$, $d_W\approx1.445482$). In terms of the extremal disc whose $18$-gon is the fundamental domain of Figure~\ref{fig:X6domain} (the uniformization by $K_{18}$, with an edge midpoint on the negative imaginary axis), the six Weierstrass points lie inside that disc, at distance $R_9$ from its centre in the directions $60^\circ,180^\circ,-60^\circ$ and at distance $d_W$ in the directions $0^\circ,120^\circ,-120^\circ$.
\end{enumerate}
\end{theorem}

\begin{proof}
The three groups are shown to uniformize the same surface by the criterion of Section~\ref{sec:hidden}: the groups $K_j$ attached to the hidden packings of each of them are $\Dpm$-conjugate to representatives of the other classes. In particular the group $K_{18}<\Dl(2,3,18)$ of the $1$-extremal surface is contained in five conjugates of $\Dl(2,3,9)$ (here the equality $L_9=L_{18}$ is used to carry out the verification of Step 5 in a common field), four of them conjugate to $K_9'$ and one to $K_9$.

(i) For $K=K_9$ the normalizer $N=N_{\Dpm}(K)$ has index $3$ in $\Dpm(2,3,9)$, so $|N/K|=24$, and none of the four hidden packings has $K_j$ conjugate to $K$; hence no further isometries exist. The induced conformal group $N^+/K$, $N^+=N\cap\Dl$, has twelve elements with orders $1,2^7,3^2,6^2$, hence is dihedral of order $12$. For the other two uniformizing groups ($K_9'$ and $K_{18}$) only a group of order $6$ is induced and the remaining $18$ isometries are hidden; the orbit structure of the packings follows from the list of those $j$ with $K_j\sim K$.

(ii) The quotient $X_6/\Aut^+(X_6)=\D/N^+$ is the orbifold of the subgroup $N^+$ of index $3$ in $\Dl(2,3,9)$. Its coset action is $x\mapsto\mathrm{id}$, $y\mapsto(1\,2\,3)$, $xy\mapsto(1\,2\,3)$, so $N^+$ has signature $(0;2,2,2,3)$ and $\D/N^+\to\D/\Dl=\mathbb P^1$ is the map $t=w^3$ in suitable coordinates, totally ramified over the points of orders $3$ and $9$; the three cone points of order $2$ are the preimages of the point of order $2$, i.e.\ the cube roots of unity, and the cone point of order $3$ is $w=\infty$. This configuration $\{1,\omega,\omega^2,\infty\}$ is equianharmonic (its cross ratio is a primitive sixth root of unity): it is invariant under the deck transformation $w\mapsto\omega w$. On the other hand, every genus two curve with $D_6\subset\Aut^+$ has an equation $y^2=x^6+ax^3+1$ with $\psi\colon x\mapsto\omega x$, $\xi\colon x\mapsto 1/x$ and $J\colon y\mapsto -y$ \cite{Bolza_1887, Broughton_1991}; the quotient map by $D_6$ is $w'=x^3+x^{-3}$ and the cone points are $w'=2,-2$ (fixed points of the non-central involutions), $w'=-a$ (Weierstrass points) and $w'=\infty$ (fixed points of $\psi$, order $3$). Since the cone point of order $3$ is at $\infty$, the configuration is equianharmonic if and only if the three cone points of order $2$ form an equilateral triangle (a M\"obius transformation of order $3$ fixing $\infty$ is a rotation by $2\pi/3$), and the triangle $\{2,-2,-a\}$ is equilateral if and only if $a=\pm2\sqrt3\,i$; the two signs give isomorphic curves ($x\mapsto -x$). The roots of $r^2+2\sqrt3 i\,r+1$ are $r=(2-\sqrt3)i$ and $r=-(2+\sqrt3)i$, and the substitution $x=\lambda x'$ with $\lambda^3=(2-\sqrt3)i$ gives the second equation. Since $X_6$ is real, it also has a model with real coefficients: the substitution $x=\mu x'$, $y=iy'$ with $\mu^3=i$ turns $y^2=x^6+2\sqrt3\,i\,x^3+1$ into $y'^2=x'^6+2\sqrt3\,x'^3-1$.

The same value of $a$ is obtained independently from the uniformization by $K_{18}$, for which only the cyclic group $C_3=\langle\psi\rangle$ is induced: the quotient $X_6/C_3=\D/N_{\Dl}(K_{18})$ is a sphere with four cone points of order $3$, its Belyi map to $\D/\Dl(2,3,18)$ has degree $6$ with cycle types $(2,2,2)$, $(3,1,1,1)$, $(6)$, which forces $t=(u^3+1)^2$, with cone points at $u=\infty$ and $u^3=-2$, again an equianharmonic configuration; in the family $y^2=x^6+ax^3+1$ the four fixed points of $\psi$ have cross ratio $(a-2)/(a+2)$ on $X_6/C_3$, which is a primitive sixth root of unity exactly when $a=\pm2\sqrt3 i$.

(iii) For $K=K_9$ the element $ab$ normalizes $K$ and the induced involution has six fixed points, all of them midpoints of edges, so it is the hyperelliptic involution; the six other conformal involutions have two fixed points each. These fixed points are determined exactly: since for $K_9$ all the isometries are induced, i.e.\ the normalizer of $K$ in $\Isom(\D)$ is contained in $\Dpm(2,3,9)$, a fixed point of an involution of $X_6$ lifts to the fixed point of an element of order $2$ of $\Dl(2,3,9)$, that is, to the midpoint $g(m_0)$ of an edge, $m_0$ being the vertex of $T_0$ with angle $\pi/2$; and $g(m_0)$ projects to a fixed point of the involution induced by $ab$ if and only if $K\,ab\,g\in\{Kg,\,Kg\,ab\}$, a condition on the permutation representation. The distance $d_W$ is computed in the triangle formed by $O$, the centre $O'$ of the $9$-gon adjacent to $P$ across an edge ($\rho(O,O')=2R_9$) and the midpoint of an edge of that $9$-gon adjacent to the common edge, which subtends the angle $2\pi/9$ at $O'$. The positions for the uniformization by $K_{18}$ are the images of these points under the isometry conjugating the two groups, which is one of the exact matrices $\alpha$ of Step~5. The three Weierstrass points at distance $R_9$ from $O$ in this second picture are consistent with \cite{Girondo_Gonzalez-Diez_2002_2}: there the lift of $J$ maps $o$ to the centre $w$ of another extremal disc, a vertex of the $18$-gon, so the midpoint of the segment $[o,w]$ is a Weierstrass point (the solution of $J(z)=z$ in Proposition~1), at distance $\rc(18)/2=R_9$ from $o$, since $\cosh\rc(18)=\cot\frac{\pi}{18}\cot\frac{\pi}{3}=\cosh 2R_9$.
\end{proof}

\begin{figure}[h]
\centering
\includegraphics[width=0.8\textwidth]{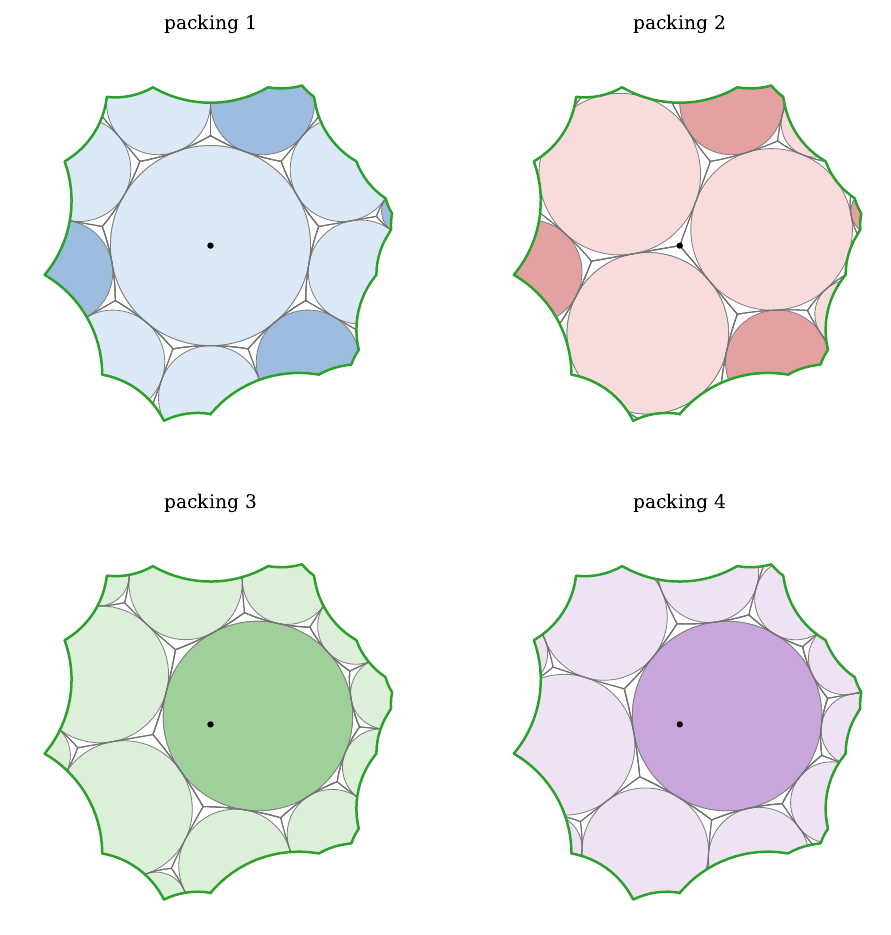}
\caption{The four extremal $4$-packings of the surface $X_3$ of Remark~\ref{rem:21} (\texttt{fuchs\_N9\#21}), drawn in one and the same fundamental domain: the Dirichlet domain (green) centred at a disc of the packing that comes with the uniformization (packing 1), whose centre is marked in the four pictures. Each picture shows the tessellation by regular $9$-gons of one packing and its discs, cut along the sides of the domain. The darker discs are those shared with another packing: packings 1 and 2 share one disc, and so do packings 3 and 4; the shared disc of the second pair is centred at distance $\rc(9)/2$ from the marked point.}\label{fig:X3}
\end{figure}

\begin{remark}\label{rem:21}
The only other surface with more than two packings is the $4$-extremal surface $X_3$ ($N=9$; \texttt{fuchs\_N9\#21}) with four packings; it is also $1$-extremal with two extremal discs, and as such it is the surface $X_3$ of \cite{Girondo_Gonzalez-Diez_2002_2} (Theorem~\ref{thm:or}(v)). Its full conformal automorphism group is $V_4$ and its isometry group has order $8$, the only isometry induced by the triangle group being an anticonformal involution; the four packings form a single orbit, with stabilizer the induced anticonformal involution; they split into two pairs, exchanged by the isometries, and the two packings of each pair share one disc (for one of the pairs its centre lies at distance $\rc(9)/2$, half the circumradius of the $9$-gons, from the centre of the fundamental polygon, and is a fixed point of a conformal involution). The four packings are shown in Figure~\ref{fig:X3}. The hyperelliptic involution is hidden; its six fixed points are fixed points of explicit elliptic elements of $\PGL_2(L_9)$, and lie at distances approximately $1.034$, $1.283$, $1.754$, $2.066$, $2.559$ and $2.709$ from the centre of the fundamental polygon. The curve is of the form $y^2=x^6+bx^4+cx^2+1$; since $V_4$ is not induced by the triangle group, the orbifold $X_3/V_4$, a sphere with five cone points of order $2$, does not cover $\D/\Dl(2,3,9)$ and the determination of $b,c$ requires the numerical uniformization of this orbifold, which we have not carried out.
\end{remark}

\begin{remark}[Comparison with \cite{Girondo_Gonzalez-Diez_2002_2}]\label{rem:ggd}
The case $k=1$, $N=18$ was treated in \cite{Girondo_Gonzalez-Diez_2002_2} by hand, starting from the nine side pairings of the $18$-gon found by Fricke and Klein. Our results agree with that paper in every detail. The eight $\Dpm$-classes of Table~\ref{ta:all_info} correspond to the surfaces $X_1,\dots,X_8$ of \cite{Girondo_Gonzalez-Diez_2002_2} (the identification of the classes with the polygons $P_1,\dots,P_8$ of Fricke and Klein is made by comparing the side pairing patterns up to rotations and reflections of the $18$-gon) and the ninth $\Dl$-class to the mirror image $\overline{X_8}$ of $X_8$, which is the only one of the eight surfaces without anticonformal automorphisms (Theorem~\ref{thm:or}(ii)). The number of extremal discs ($1$ for $X_2$, $4$ for $X_6$, $2$ otherwise), the full automorphism groups ($\Z_2$ for $X_2,X_4,X_5,X_7$; $\Z_2\times\Z_2$ for $X_1,X_3,X_8,\overline{X_8}$; $D_6$ for $X_6$) and the positions of the hidden centres (the midpoint $p$ of a side for $X_1,X_3,X_8,\overline{X_8}$; the points $q$ at distance approximately $0.777654$, $1.28738$ and $1.21231$ from the centre for $X_4$, $X_5$ and $X_7$; the three vertices $u,v,w$ for $X_6$) are all reproduced by the algorithm, including the automorphisms that are not induced by the triangle group, which in \cite{Girondo_Gonzalez-Diez_2002_2} were found by a case-by-case geometric analysis. Lemma 1 of \cite{Girondo_Gonzalez-Diez_2002_2} shows that $X_6$ is a curve $y^2=(x^3-1)(x^3-\lambda^6)$ for some $\lambda$ with $\lambda^6\neq 1$, and Theorem~\ref{thm:X6} determines the parameter: $\lambda^6=-(2+\sqrt3)^2=-(7+4\sqrt3)$. \end{remark}

\begin{remark}[Comparison with \cite{Girondo_Torres_Wolfart_2012}]\label{rem:gtw}
The surface $X_6$ appears for a third time in Section 7 of \cite{Girondo_Torres_Wolfart_2012}, as a genus two surface $S\backslash\mathbb H$ with $S\lhd\Dl(\pi)<\Dl(3,3,9)<\Dl(2,3,9)$, where $\pi$ is the ramified prime of $\Q(\cos\frac{\pi}{9})$ dividing $3$ and $\Dl(\pi)$ is the principal congruence subgroup, of signature $(0;3,3,3,3)$. Uniform dessins of type $(2,3,9)$ on a surface of genus $2$ are the same thing as extremal $4$-packings, and the description given there --- four $(2,3,9)$-dessins equivalent under the automorphism group $N(S)/S\cong D_3\times C_2$ plus one stabilized by it, corresponding to the middle vertex of the subtree of local maximal orders containing $S$ --- is exactly the orbit structure of Theorem~\ref{thm:X6}(i); the arithmetic of the Eichler orders of $\Dl(2,3,9)$ explains \emph{a priori} why $X_6$ is the surface with the largest number of packings in our list. The trigonal equation $y^3=(x-1)(x^3-1)$ quoted there from \cite{Singerman_Syddall_2003} is consistent with Theorem~\ref{thm:X6}(ii): the quotient by the cyclic group of order $3$ is branched over $1,\omega,\omega^2,\infty$, an equianharmonic configuration, and in the family $y^2=x^6+ax^3+1$ the four fixed points of $\psi$ have cross ratio $(a-2)/(a+2)$ on $X/\langle\psi\rangle$, which is a primitive sixth root of unity exactly for $a=\pm2\sqrt3\,i$. On the other hand, the hyperelliptic model $y^2=x^6+8x^3+4$ given in Remark 5 of \cite{Girondo_Torres_Wolfart_2012} seems to contain a misprint: after scaling it is the member $a=4$ of the family, whose cross ratio $1/3$ is not equianharmonic (its Igusa--Clebsch invariants differ from those of $y^2=x^6+2\sqrt3\,i\,x^3+1$, for instance $A^5/D=-49152$ versus $A^5/D\approx 990114.67$); the correct equation is the one of Theorem~\ref{thm:X6}(ii).
\end{remark}

\begin{remark}\label{rem:autgroups}
Since every isometry of a surface is obtained explicitly, as a matrix of $\PGL_2(\R)$ normalizing $K$, the isomorphism types of $\Isom(X)$ and, in the orientable case, of $\Aut^+(X)$ can be determined for all the surfaces: the multiplication table of these matrices modulo $K$ defines the group, which is then identified with GAP. For the orientable surfaces of genus $2$ the group $\Aut^+(X)$ is $\Z_2$, $\Z_2\times\Z_2$, $D_4$, $D_6$, $\Z_{10}$ or $\mathrm{GL}(2,3)$, the last two only for the curves $y^2=x^5-1$ and $y^2=x^5-x$ of Remark~\ref{rem:bolza} (as expected, since these are their well known automorphism groups). For the non-orientable surfaces of genus $3$ the isometry group is $\Z_2$, $\Z_2\times\Z_2$, $D_4$ or $D_6$. The groups of every surface are listed in the data files (Section~\ref{sec:data}).
\end{remark}

\begin{figure}[h]
\centering
\includegraphics[width=0.8\textwidth]{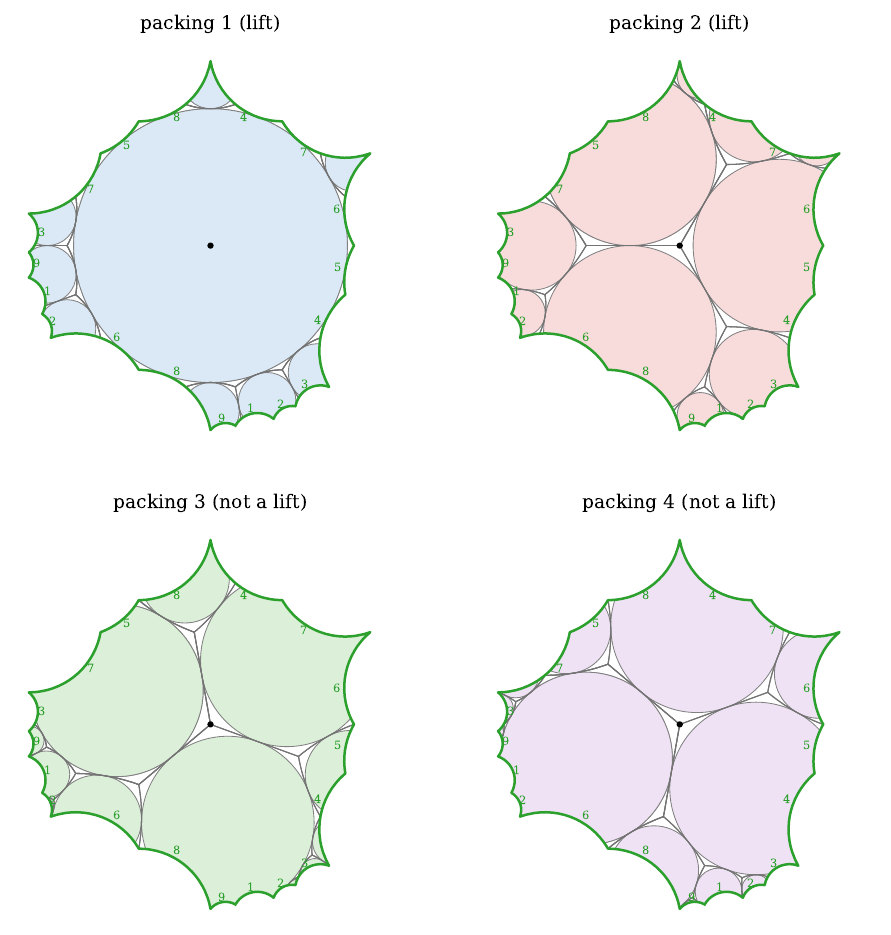}
\caption{The four extremal $2$-packings of the genus three surface $Y$ of Remark~\ref{rem:genus4}, which is the canonical double cover of two non-isometric non-orientable surfaces of genus $4$, drawn in one and the same Dirichlet domain (green), centred at a disc of packing 1. Packings 1 and 2 are the lifts of the two extremal discs of one of them; packings 3 and 4, exchanged by its covering involution, are the lifts of the discs of the other.}\label{fig:genus4}
\end{figure}

\begin{remark}[Double covers in genus $4$]\label{rem:genus4}
As recalled after Theorem~\ref{thm:cover}, in general a double cover has at least as many extremal packings as the non-orientable surface: an extremal $2k$-packing of $X^+$ descends to $X$ if and only if it is invariant under the covering involution $\tau$, and those that are not invariant come in pairs. The equality may fail outside genus $3$. We have run the algorithm on the $144$ non-orientable $1$-extremal surfaces of genus $4$ classified in \cite{Nakamura_2009} and on their double covers, which are $2$-extremal surfaces of genus $3$. In five cases, exactly those whose second extremal disc is centred at a vertex of the $18$-gon, the double cover has four extremal $2$-packings instead of two: the lifts of the two discs of $X$, and two packings exchanged by $\tau$. These are lifts with respect to another covering involution. In three of the cases this involution is conjugate to $\tau$ by an isometry of $X^+$, and the four packings form a single orbit. In the other two, the two non-orientable surfaces have the same double cover $Y$, whose four packings form two orbits, the lifts of the discs of each quotient (Figure~\ref{fig:genus4}). Besides $X_6$, $Y$ is therefore an example of a surface uniformized by two subgroups of the same triangle group that are not conjugate in it (Lemma~\ref{lem:orbits}). The computation is done by the script \texttt{genus4\_covers.py} of the repository (Section~\ref{sec:data}), whose output lists these surfaces.
\end{remark}

\section{Data}\label{sec:data}
\begin{figure}[h]
\centering
\includegraphics[width=0.9\textwidth]{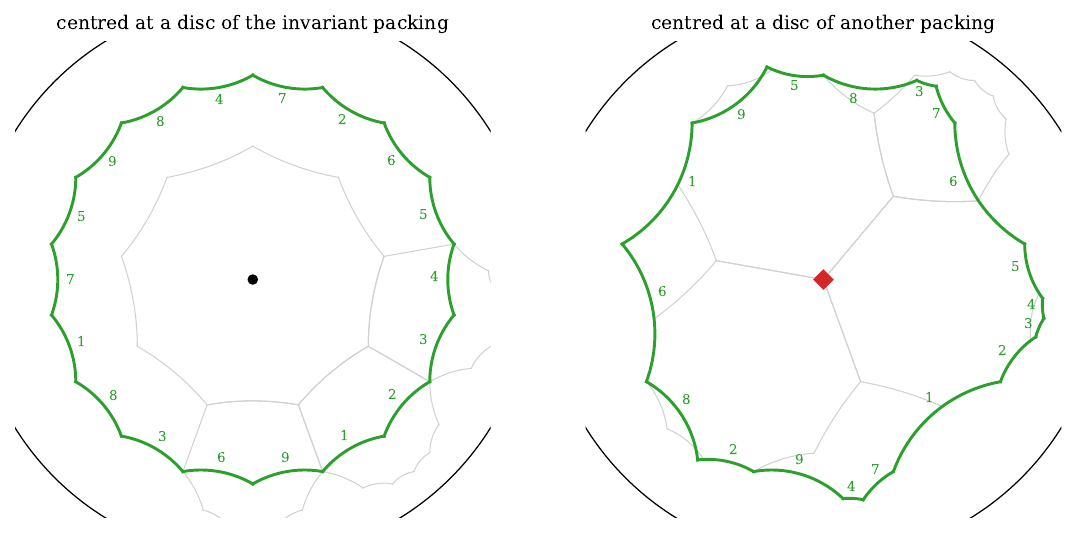}
\caption{Two Dirichlet domains of $X_6$ (green; sides with the same number are identified), uniformized by the group $K_9$ of Theorem~\ref{thm:X6} (so that the $9$-gons of the standard tessellation are those of the invariant $4$-packing), with the four $9$-gons of the fundamental domain of Step~1 in grey. Left: centred at a disc of the $4$-packing invariant under the isometry group, whose centre is also the centre of one of the four extremal discs of the $1$-extremal surface $X_6$; the domain is the regular $18$-gon. Right: centred at a disc of one of the other four $4$-packings. Each picture is centred at the centre of its domain.}\label{fig:dirichlet}
\end{figure}

The data and the programs described in this section allow the reader both to access all the information of the catalogue and to reproduce it completely; every computation can be repeated from the enumeration of the groups, and a script compares the result with the distributed data file by file. For each of the $423$ surfaces (the $373$ of Table~\ref{ta:all_info} and the $50$ double covers) the following data are stored: the permutation representation of $\Dpm$ on the cosets of $K$ and generators of $K$ as given by GAP; the side pairings of the fundamental domain as words in $a,b,c$; for each extremal packing, the coordinates of a centre in the disc model, the exact matrix $\alpha\in\PGL_2(L_N)$ and the words $w_s$ of Proposition~\ref{prop:cert}; the automorphisms induced by $\Dpm$ (coset, parity, order, existence of fixed points); the canonical form of the permutation representation, which identifies the $\Dpm$-conjugacy class; the list of hidden isomorphisms; and the isomorphism types of $\Isom(X)$ and $\Aut^+(X)$ (Remark~\ref{rem:autgroups}). Programs and data are available at \url{https://github.com/ernesto-girondo/girondo_reyes_lowgenus_packings} and are archived at Zenodo \cite{Girondo_Reyes_data}. All the software used (GAP \cite{GAP4}, PARI/GP \cite{PARI2}, and Python with the libraries \texttt{mpmath} \cite{mpmath} and \texttt{cypari2}) is free and open source, so that every computation can be repeated without any commercial software; the file \texttt{VERIFYING.md} there explains, statement by statement, how to check the results of this paper. Besides the data files, the programs draw, for any of the $373$ groups, three pictures, each one accompanied by a text file with the approximate coordinates (in the disc model, together with the hyperbolic distance to $O$ and the angle) of all the points shown. The first picture shows the fundamental domain of Step~1, with its $k$ polygons labelled and the edges identified by each side pairing numbered (orientation reversing pairings being marked); the centres of the discs of all the extremal packings; and the fixed points and mirrors of all the isometries of the surface, induced or hidden, the fixed points being marked according to the order of their stabilizers and the Weierstrass points being singled out (for every orientable surface of genus $2$ the text file gives the six Weierstrass points, as the fixed points of the hyperelliptic involution, whether this involution is induced by the triangle group or hidden). The text file lists, in addition, the side pairings as words in $a,b,c$, the isometries with their orders and numbers of fixed points, and the discs shared by two different packings. The second picture shows, for each packing, the Dirichlet domain of $K$ centred at one of its discs, with its side pairings. The third one shows Step~3 of the search, as Figure~\ref{fig:edcycles} does for $X_6$: the ED-cycles of the two side pairings used, all their intersection points in the search disc and the candidates, each candidate marked as the centre of a disc of a packing or as a point shown not to be a centre; its text file lists every admissible distance up to $D$ and every intersection point, with the ED-cycles on which it lies and the reason why it is kept or discarded, so that the search can be followed step by step for any surface. It is computed from the orbit of the centre, which is obtained exactly from the permutation representation (a centre $g(O)$ of the tessellation lies in $K\cdot O$ if and only if $gr\in K$ for some $r$ in the stabilizer $\langle a,c\rangle$ of $O$), as an intersection of half planes, using all the orbit points up to twice the distance from the centre to the farthest vertex, so that no side is missing; every side is paired, every vertex cycle has angle sum $2\pi$ and the area is $2\pi|\chi|$, so that by Poincar\'e's theorem the side pairings generate $K$. For instance, the Dirichlet domain of $X_6$ centred at a point that is the centre of a disc of the invariant $4$-packing and also of one of the four extremal discs (Theorem~\ref{thm:X6}) is the regular $18$-gon with all its vertices at distance $\rc(18)$ from its centre, while the one centred at a disc of one of the other four packings is an $18$-gon with vertices at three different distances from its centre (Figure~\ref{fig:dirichlet}). Over all the surfaces, the Dirichlet domain centred at $O$ has between $8$ and $18$ sides for the orientable surfaces and between $8$ and $12$ for the non-orientable ones; the classical regular octagon of Bolza's surface is one of the $17$ octagons with a single vertex cycle. A program draws all the pictures and checks, surface by surface, the numbers of packings, the orders of the induced and full isometry groups, the integrality of $\chi(X/G)$ computed from the fixed points of the isometries, and the conditions above on the Dirichlet domains, and that every candidate of Step~3 is either the centre of a disc of a packing or shown not to be one.

 Surfaces are referred to by their family and their position in the GAP enumeration, as explained in Section~\ref{sec:classes} (for instance \texttt{fuchs\_N9\#21} is the surface $X_3$ of Remark~\ref{rem:21}), which is the numbering used throughout the paper. In the data files the index of a group is always its index in $\Dpm$, which for an orientable surface is twice the index $d$ of Section~\ref{sec:groups}. This order is the one given by GAP 4.12.1 and 4.15.1, and another version of GAP might list the classes in a different order. For this reason the data files also contain a canonical form of each class, which identifies it up to conjugation in $\Dpm$ whatever the order, and a short program identifies any surface given by its monodromy (the permutations of $a,b,c$, or the monodromy of its dessin of type $(2,3,N)$) in the catalogue; a computation that produces the classes in another order can therefore be compared with ours class by class. The existence part of every statement can be checked without running the search: a short independent script checks, in exact arithmetic over $L_N$, the identities that confirm the $646$ packings (for each extremal packing, $s\alpha=\pm\alpha w_s$ for every side pairing $s$, together with the facts that the side pairings lie in $K$ and generate it, and that the packings are pairwise different). The completeness part (that there are no further packings) rests on the search of Section~\ref{sec:algorithm}, and the reader who wishes to check it can repeat the whole computation, from the enumeration of the groups to the tables of Section~\ref{sec:results}, in less than an hour on a standard computer.

\subsection*{Acknowledgements}
This work was partially supported by the Severo Ochoa Programme for Centres of Excellence in R\&D (CEX2023-001347-S), funded by MCIN/\allowbreak AEI/\allowbreak 10.13039/\allowbreak 501100011033, Spain, and by the project ``Extremos de estructuras geom\'etricas en superficies hiperb\'olicas de tipo finito'' (UAM Avanza 2026, ref.\ 2026/00263/001), funded by the Universidad Aut\'onoma de Madrid.

\section*{Use of AI tools}
The computational part of this work was prepared with the assistance of Claude (Anthropic), an AI model (Claude Opus 5.5), used between October and November 2026. The authors provided the method (\cite[Chapter~7]{Reyes_thesis}) together with earlier versions of the programs, which already carried out the computations surface by surface. Under the direction of the authors, Claude organized, cleaned and restructured these programs into the single, documented pipeline described in Sections~\ref{sec:algorithm} and~\ref{sec:data}, ran the computations on all the surfaces, organized the resulting data, wrote, at the request of the authors and in order to give more complete information, the auxiliary programs (pictures of the surfaces, Dirichlet domains, identification of the automorphism groups and of a surface given by its monodromy), and prepared the documentation that allows an interested reader to reproduce and check every computation; it also proposed revisions of the text. The authors checked the results against the literature (Section~\ref{sec:algorithm}) and through the independent exact verification described in Section~\ref{sec:data}, and take full responsibility for the content of the paper, including its references.

\appendix
\section{Canonical double covers}\label{app:covers}
The following table lists, for each non-orientable $k$-extremal surface of genus $3$ (numbered as in the GAP enumeration of NEC subgroups), the Fuchsian class of its canonical double cover (numbered as in the GAP enumeration of subgroups of $\Dl(2,3,N)$), the number $p$ of extremal packings and the order $A$ of the group of isometries induced by the triangle group. The numbering of the classes may differ from the one used in \cite{Reyes_thesis}: for $N=9$ the NEC class \texttt{nec\#2} is the group $K_{(9,-,2)}$ of Section~7.1 of \cite{Reyes_thesis} (the permutations printed there define the same class, with a different numbering of the cosets), but the Fuchsian classes were enumerated differently; for instance, the double cover of \texttt{nec\#2} for $N=9$, called $K_{(9,+,13)}$ in Section~7.1 of \cite{Reyes_thesis}, is \texttt{fuchs\#20} here.

\begin{center}\small
\begin{longtable}{llcc@{\qquad}llcc}
\toprule
NEC & cover & $p$ & $A$ (NEC/cover) & NEC & cover & $p$ & $A$ (NEC/cover)\\
\midrule
\endhead
\multicolumn{8}{l}{$N=12$, $k=1$ (covers are $2$-extremal)}\\
nec\#1 & fuchs\#5 & 2 & 2/4 & nec\#7 & fuchs\#2 & 1 & 2/4\\
nec\#2 & fuchs\#9 & 2 & 1/2 & nec\#8 & fuchs\#4 & 2 & 2/4\\
nec\#3 & fuchs\#8 & 2 & 2/12 & nec\#9 & fuchs\#13 & 2 & 2/4\\
nec\#4 & fuchs\#11 & 2 & 6/12 & nec\#10 & fuchs\#14 & 2 & 2/4\\
nec\#5 & fuchs\#1 & 2 & 2/4 & nec\#11 & fuchs\#24 & 2 & 6/12\\
nec\#6 & fuchs\#21 & 2 & 2/4 & & & & \\
\midrule
\multicolumn{8}{l}{$N=9$, $k=2$ (covers are $4$-extremal)}\\
nec\#1 & fuchs\#12 & 1 & 2/4 & nec\#3 & fuchs\#17 & 2 & 2/4\\
nec\#2 & fuchs\#20 & 2 & 2/4 & nec\#4 & fuchs\#16 & 2 & 1/2\\
\midrule
\multicolumn{8}{l}{$N=8$, $k=3$ (covers are $6$-extremal)}\\
nec\#1 & fuchs\#18 & 2 & 2/4 & nec\#13 & fuchs\#9 & 1 & 2/4\\
nec\#2 & fuchs\#1 & 2 & 2/4 & nec\#14 & fuchs\#5 & 2 & 2/4\\
nec\#3 & fuchs\#16 & 2 & 6/12 & nec\#15 & fuchs\#3 & 1 & 2/4\\
nec\#4 & fuchs\#17 & 2 & 2/12 & nec\#16 & fuchs\#12 & 2 & 2/4\\
nec\#5 & fuchs\#27 & 2 & 2/4 & nec\#17 & fuchs\#75 & 2 & 2/4\\
nec\#6 & fuchs\#26 & 2 & 1/2 & nec\#18 & fuchs\#38 & 1 & 4/8\\
nec\#7 & fuchs\#28 & 1 & 2/4 & nec\#19 & fuchs\#42 & 1 & 4/8\\
nec\#8 & fuchs\#59 & 2 & 2/4 & nec\#20 & fuchs\#55 & 1 & 4/8\\
nec\#9 & fuchs\#31 & 2 & 2/4 & nec\#21 & fuchs\#51 & 2 & 6/12\\
nec\#10 & fuchs\#39 & 1 & 2/4 & nec\#22 & fuchs\#64 & 1 & 8/16\\
nec\#11 & fuchs\#58 & 2 & 6/12 & nec\#23 & fuchs\#76 & 2 & 6/12\\
nec\#12 & fuchs\#61 & 2 & 2/12 & & & & \\
\midrule
\multicolumn{8}{l}{$N=7$, $k=6$ (covers are $12$-extremal)}\\
nec\#1 & fuchs\#46 & 2 & 1/2 & nec\#7 & fuchs\#55 & 2 & 1/2\\
nec\#2 & fuchs\#1 & 2 & 2/4 & nec\#8 & fuchs\#50 & 1 & 2/4\\
nec\#3 & fuchs\#41 & 2 & 6/12 & nec\#9 & fuchs\#2 & 2 & 2/4\\
nec\#4 & fuchs\#42 & 2 & 2/12 & nec\#10 & fuchs\#140 & 2 & 2/4\\
nec\#5 & fuchs\#44 & 1 & 2/4 & nec\#11 & fuchs\#132 & 1 & 2/4\\
nec\#6 & fuchs\#54 & 2 & 2/4 & nec\#12 & fuchs\#128 & 2 & 6/12\\
\bottomrule
\end{longtable}
\end{center}

\bibliographystyle{alpha}
\bibliography{references}

\end{document}